\documentclass[letterpaper,12pt]{article}
\usepackage{times}
\usepackage{empheq}
\usepackage{amsmath, amssymb}
\usepackage{amsthm}
\numberwithin{equation}{section}
\usepackage{bbm}
\usepackage[T1]{fontenc}
\usepackage[english]{babel}
\usepackage[utf8]{inputenc}
\usepackage{hyperref}
\usepackage{mathrsfs}
\usepackage{mathabx}
\usepackage{color}
\usepackage{anysize}
\usepackage{apacite}
\usepackage{natbib}
\marginsize{2cm}{2cm}{2cm}{2cm}
\usepackage[toc,page]{appendix}
\usepackage{url}
\usepackage{mathtools}
\usepackage[font={small,sc}]{caption}
\usepackage{ragged2e}
\usepackage[shortlabels]{enumitem}

\newtheorem{Theo}{Theorem}[section]
\newtheorem{Def}{Definition}[section]
\newtheorem{Prop}{Proposition}[section]
\newtheorem{Corol}{Corollary}[section]
\newtheorem{Lemma}{Lemma}[section]
\newtheorem{Rem}{Remark}[section]
\theoremstyle{definition}
\newtheorem{Exa}{Example}[section]

\AtBeginDocument{}

\DeclareMathOperator\sgn{sgn}

\DeclareMathOperator\supp{supp}
\DeclareMathOperator\diam{diam}

\newcommand{\Rd}{\mathbb{R}^{d}}

\newcommand{\RdxRd}{\Rd\times\Rd}

\newcommand{\BorelRd}{\mathcal{B}(\mathbb{R}^{d})}

\newcommand{\BoundedBorelRd}{\mathcal{B}_{B}(\mathbb{R}^{d})}

\newcommand{\Cov}{\mathbb{C}ov}
\newcommand{\Var}{\mathbb{V}ar}
\newcommand{\OmegaAP}{(\Omega , \mathcal{A} , \mathbb{P})}
\newcommand{\LtwoOmega}{L^{2}(\Omega , \mathcal{A} , \mathbb{P})}

\newcommand{\R}{\mathbb{R}}

\newcommand{\DRd}{\mathscr{D}(\Rd)}

\newcommand{\DistributionsRd}{\mathscr{D}'(\Rd)}

\newcommand{\BoundedBorelRp}{\mathcal{B}_{B}(\R^{p})}

\newcommand{\RiemannPartitionTag}[3]{\left( (#1_{j}^{n})_{j \in #3_{n}} , (#2_{j}^{n})_{j \in #3_{n}} \right)_{n \in \mathbb{N}}}
\newcommand{\RiemannSystem}[3]{\left( (#1_{j}^{n})_{j \in #3_{n}} , (#2_{j}^{n})_{j \in #3_{n}} \right)_{n \in \mathbb{N}}}
\newcommand{\Borel}[1]{\mathcal{B}(#1)}
\newcommand{\BoundedBorel}[1]{\mathcal{B}_{B}(#1)}

\newcommand{\LpOmegaAP}[1]{L^{#1}(\Omega , \mathcal{A} , \mathbb{P})}
\newcommand{\LpOmega}[1]{L^{#1}(\Omega , \mathcal{A} , \mathbb{P})}

\begin{document}

\begin{Large}
\begin{center}
\textbf{Function-measure kernels, self-integrability and uniquely-defined stochastic integrals}
\end{center}
\end{Large}
\begin{center}
\large
\textsc{Ricardo Carrizo Vergara} \\
\normalsize
\textit{racarriz@uc.cl}
\end{center} 
\begin{center}
March 2023
\end{center}

\begin{center}
\textbf{Abstract}
\end{center}

\small
\justify
In this work we study the self-integral of a function-measure kernel and its importance on stochastic integration. A continuous-function measure kernel $K$ over $D \subset \Rd$ is a function of two variables  which acts as a continuous function in the first variable and as a real Radon measure in the second. Some analytical properties of such kernels are studied, particularly in the case of cross-positive-definite type kernels. The self-integral of $K$ over a bounded set $D$ is the ``integral of $K$ with respect to itself''. It is defined using Riemann sums and denoted $\int_{D}K(x,dx)$. Some examples where such notion is well-defined are presented. This concept turns out to be crucial for unique-definiteness of stochastic integrals, that is, when the integral is independent of the way of approaching it. If $K$ is the cross-covariance kernel between a mean-square continuous stochastic process $Z$ and a random measure with measure covariance structure $M$, $\int_{D}K(x,dx)$ is the expectation of the stochastic integral $\int_{D} ZdM$ when both are uniquely-defined. It is also proven that when $Z$ and $M$ are jointly Gaussian, self-integrability properties on $K$ are necessary and sufficient to guarantee the unique-definiteness of $\int_{D}ZdM$. Results on integrations over subsets, as well as potential $\sigma$-additive structures are obtained. Three applications of these results are proposed, involving tensor products of Gaussian random measures, the study of a uniquely-defined stochastic integral with respect to fractional Brownian motion with Hurst index $H > \frac{1}{2}$, and the non-uniquely-defined stochastic integrals with respect to orthogonal random measures. The studied stochastic integrals are defined with great generality on the integrand and the integrator and without use of martingale-type conditions, providing a potential filtration-free approach to stochastic calculus grounded on covariance structures.
\bigskip


\normalsize

\section{Introduction}
\label{Sec:Introduction}

This work is focused on a concept which, to the best of our knowledge, has not received much attention in the literature, neither as an interesting mathematical object it its own right, nor as a potentially important tool in applied mathematics. This concept is the \textit{self-integral of a function-measure kernel}, and its main application is in stochastic integration. A function-measure kernel is a structure defined over the cross product of two measurable spaces which acts as a function in one component and as a measure in the other.  The \textit{self-integral} is a sort of \textit{integral of the kernel with respect to itself}, that is, of its \textit{first component} with respect to its \textit{second component}. Such a notion may seem peculiar at a first sight but it is simple to conceive examples where it is well-defined. The purpose of this work is triple: to study some basic analytical properties of function-measure kernels, to define and study the concept of self-integrability, and to show the importance of such concept in stochastic integration, where it is found to be a key concept regarding uniquely-defined stochastic integrals.

We consider a function $K : D \times \mathscr{E} \to \R $, with $D$ some measurable space and $\mathscr{E}$ some $\sigma$-algebra (or $\delta$-ring), for which $K(\cdot , A)$ is a function and $K(x, \cdot)$ is a measure. Such kind of kernel is widely used in the literature, and maybe the most important example is a transition or conditional probability kernel, in which case $K(x , \cdot)$ is a probability measure and $K(\cdot , A)$ is a positive measurable function. In this work though the kernel $K$ will be used for other purposes. First of all, we consider different regularity properties: $K(\cdot , A)$ will be taken to be a continuous function over $D$ and $K(x , \cdot)$ will be a real Radon measure. We shall focus on the case $D \subset \Rd$ and $\mathscr{E} = \Borel{\R^{p}}$, though many of our results can be generalized to the classical settings where Radon measures are used (replace $D$ or $\R^{p}$ by any locally compact Hausdorff topological space which can be written as a countable union of compact sub-spaces). The main role of such kernel $K$ in a stochastic context is to describe \textit{cross-covariance structures}. This implies the use of kernels $K$ which have an extra cross-positive-definiteness condition which allows them to describe the covariance between a mean-square continuous stochastic process and a second-order random measure. 

The notion of self-integral can be conceived when  $\mathscr{E} = \Borel{D}$, hence the left and the right components of $K$ act over the same measurable space. The objective is to define the ``integral of $K$ with respect to itself'', which in this work is noted
\begin{equation}
\label{Eq:SelfIntegralIntro}
\int_{D}K(x,dx).
\end{equation}
The author is not aware of any literature which has tried to formalise such notion. Nonetheless, some readers specialised in analysis may recognize that the self-integral is a \textit{trace operator} \citep[Section 1.3]{ryan2002introduction}. Indeed, members of the tensor product space $C(D)\otimes \mathscr{M}(D)$, where $C(D)$ is the space of continuous functions over $D$ and $\mathscr{M}(D)$ is the space of real Radon measures over $D$, form a particular class of function-measure kernels where the self-integral can be defined trivially (see Section \ref{Sec:SelfIntegral}). Since for $D$ compact $C'(D) = \mathscr{M}(D)$ (Riesz representation), the self-integral over $D$ can be defined as a trace operator over $C(D)\otimes C'(D)$, and it can be eventually extended to more general cases using the completion of such tensor product with the projective topology \citep[Chapter 2]{ryan2002introduction}. This analytical approach will be studied further in a future paper where we will try to define the self-integral through series expansion methods. In this paper we focus on an apparently simpler approach: using Riemann sums to define the self-integral. That is, we do the same as with classical integrals: relying on Riemann partitions of the integration domain, tag-points to evaluate the left component of $K$, and so constructing a limit of Riemann sums. Such self-integral is not always uniquely-defined: the final limit may depend upon the selected partition or tag-points.  It is hence necessary to study the cases where such limit is uniquely-defined and so define the concept of self-integral for those.

The self-integral appeared naturally when trying to define uniquely the stochastic integral between a general mean-square continuous stochastic process and a second-order random measure (see Definition \ref{Def:RandomMeasure} for the exact notion of \textit{random measure} we use in this work). By \textit{uniquely}, we mean that the integral must not depend upon the way of approaching it. We use an approach which in our opinion is simpler than most part of the approaches already present in the literature: the definition using a limit of Riemann sums, and the study of unique-integrability conditions (independence upon partitions and tag-points) through the analysis of the \textit{cross-covariance structure} between the integrand and the integrator. Such cross-covariance structure can be described through a function-measure kernel, the same mathematical object described above. Here we show that the possibility of defining uniquely stochastic integrals is intimately related to the possibility of defining the self-integral of this cross-covariance kernel: it is indeed a necessary condition for it, since \textit{the mean of the stochastic integral is the self-integral of the cross-covariance kernel}. In addition, when the integrator and the integrand are jointly Gaussian, necessary and sufficient conditions for a uniquely-defined stochastic integral rely entirely on self-integrability properties of the cross-covariance kernel. Such conditions are exposed in Theorem \ref{Theo:IntZMGaussian}, which is the main result of this work. While, of course, the study of second-moment structures is not enough to fully characterise the laws of the objects, and hence it cannot cover the whole richness of the theory of stochastic processes, we believe that these results are powerful and simple enough to propose another way of tackling the subject, which may have potential both for generalization to more complicated cases, as well as for applications. We also remark that this approach can be easily applied for stochastic integrals over $\Rd$, and it is done regardless (in principle) of any martingale-type condition on the integrator or adaptability on the integrand.

\subsection{Motivations and some bibliographical revision}
\label{Sec:Motivations}

Let us expose some thoughts about the current developments on stochastic integrals. There is one simple but important question for which this work aims to give some enlightening: \textit{when can we integrate a stochastic process with respect to a random measure?} One could try to mimic the definition of the integral in the deterministic case, but it is widely known that in the stochastic world things are not that simple. Classical stochastic integration against (the derivative of) a Brownian motion uses intuitively constructed stochastic integrals through Riemann sums and simple-functions approximations, for which different ways of approaching the integral produce different results (classical difference between the Itô and Stratonovich integrals, see \citep[Chapter 3]{oksendal2003stochastic}). When can we thus integrate the stochastic process with respect to the random measure \textit{uniquely}? Which properties must satisfy both the process and the random measure in order to assure the limit to be well-defined and independent of the way of approaching it?

It turns out that this question seems not to have been considered in the literature. While the study of stochastic integrals has produced a succulent quantity of work in many directions for over 100 years since the first works of Wiener \citep{wiener1923differential}, the particular question of conditions for uniquely-defined stochastic integral seems not to have been taken into account. In \citep[Section 37.3]{loeve1978probability}, the author proposes a definition of a stochastic integral in a quite general setting of a second-order process as integrand and the increments of a second-order process whose covariance has bounded variation as integrator\footnote{The notion of random measure that we use in this work actually coincides with this one, although we explicitly work with the covariance measure rather than with the increments of its \textit{primitive} function.}. In such definition it is required that the limit must be independent of the selected tag-points in the Riemann-Stieltjes sum. However, the author does not provide general criteria for such condition to hold, except for the particular case when the integrator and the integrand are independent, case which we will revisit in Section \ref{Sec:FirstApproachStochIntegrals}. It is also known that when the integrand or the integrator is deterministic the stochastic integral is uniquely-defined, but this is actually a particular case of independence  (we shall revisit such constructions in Section \ref{Sec:RandomMeasures}). Besides such cases, it is not easy to find bibliographical sources with \textit{relatively simple} approaches to stochastic integration with huge generality and caring about a \textit{unique} definition. For instance, in \citep{huang1978stochastic} a construction of a stochastic integral with general Gaussian integrators and general (non-adapted) integrands is done following an approach based on tensor products and extensions on Hilbert spaces, but the approach has already been 
considered complicated \citep{borisov2006constructing}. Another attempt to properly define a particular form of stochastic integral regardless on the way of approaching it is the one present in \citep[Chapter 5]{rao2012random}. Here the author looks for a definition of a convolution between two second-order random measures, which implicitly uses a stochastic integral. Roughly speaking, the author describes the convolution between the \textit{covariance bi-measures} of both random measures as a new covariance bi-measure, and thanks to Kolmogorov Theorem (see Proposition 5.1.3 in \citep{rao2012random}), they construct a new random measure following such covariance. However, this is an abstract construction which does not link \textit{directly} the paths of the initial random measures with the paths of the final convolution. For instance, a similar and maybe more exploitable result can be obtained if we use the convolution of two independent second-order random measures, which can be defined without any problem.

On typical approaches on stochastic integration there is another point worth mentioning. Virtually every approach presents some special conditions on the integrand and/or the integrator which may seem more restrictive than what we could conceive at a first (maybe naive) sight. The integrator must satisfy a martingale-type condition, rather than just belonging to a simple  and general enough class of random measures. A semi-martingale condition on the primitive of the random measure seems to be the most general condition one may find \citep{kunita1967square,bichteler2006random,dellacherie2006survol}. On the other hand, the integrand must be adapted to the filtration induced by the random measure, rather than being just an enough regular process (mean-square continuous for instance). These conditions are usually asked for mathematical convenience (Lebesgue Dominated Theorem holds, see \citep{bichteler2006random,chong2015integrability}), hence mathematicians have built the whole functional theory over them. It is also possible that these conditions are preferred because they follow an intuition of ``not-looking-into-the-future" \citep[Chapter 3]{oksendal2003stochastic}. It seems that they are totally inspired from a purely temporal point of view, i.e. for stochastic integrals over the real line rather than the plane or other spaces. Nonetheless, stochastic integrals over the plane have been developed by applying those same principles to the spatial axes and thus doing iterated Itô integrals  \citep{ito1951multiple,cairoli1975stochastic}  .  This procedure is of course mathematically possible and legitimate but it seems ``\textit{artificial}'' since, in principle, there is no reason or human intuition for which a spatial process should follow an adaptability condition along the spatial axes\footnote{When someone requires to a stochastic processes over time to be adapted with respect to a filtration, there is an implicit idea of \textit{causality} or \textit{time arrow} behind this mathematical intuition. One must not forget that these notions are \textit{also} conceptual abstractions which we can or not make part of the mathematical models. We could in principle conceive mathematical models which do not consider them.}. We wonder hence if it is not possible to tackle the general problem of stochastic integration without using such concepts.  The adaptability is \textit{a way of describing the dependence between the process and the measure}, but maybe other simpler methods could be more exploitable in some circumstances. The already mentioned not-so-simple approaches in \citep{rao2012random,huang1978stochastic} are some examples which treat in a relatively complete manner the question of doing stochastic integration without adaptability conditions. Skorokhod integral \citep{skorokhod1976generalization}, which has been developed using Malliavin calculus \citep{nualart2006malliavin}, is also another successful approach to stochastic integration which does not consider adaptability conditions on the integrand, but it is usually restricted to the White Noise random measure (derivative of Brownian motion) for its construction. It can be extended to more general random measures \citep{mocioalca2005skorohod} but in general, it seems that it is necessary to \textit{fix the random measure} and \textit{then} construct the integration theory with respect to it. Hence, there is no a relatively free manipulation in both the integrand and the integrator at the same time. There are also other more sophisticated approaches to stochastic integration through the definition of \textit{stochastic products} between quite general stochastic structures. The stochastic integral is a sort of multiplication between a random function and a random measure, but theories which uses more general stochastic objects, such as integrators with controlled $p$-variations and random distributions have been developed during the last decades. The Wick product \citep{holden2009stochastic}, rough path theory and para-controlled distributions \citep{gubinelli2004controlling,gubinelli2015paracontrolled} and the theory of regularity structures \citep{hairer2014theory} are examples of developments which allow the definition of the product of abstract stochastic structures under suitable contexts. However, these magnificent theories are of a technicality and sophistication which could demotivate some readers from the applied world.

One final motivation for this work comes from a topic in applied mathematics: the SPDE approach in Geostatistics \citep{lindgren2011explicit,carrizov2018development}. This approach consists in the analysis of spatialised data-sets interpreted as the realisation of the solution of a Stochastic Partial Differential Equation (SPDE). We refer to the seminal paper \citep{lindgren2011explicit} as the first highly popular paper which exploits this link between SPDEs and Spatial Statistics techniques, and to \citep{lindgren2022spde} as a recent review on how the topic has evolved in the Spatial Statistics community. However, the SPDE approach is until now applied in a purely \textit{linear} setting: linear SPDEs and linear statistical methods based on covariance structures. So, there is an open door for an ambitious research program: the systematic study and exploitation of a potential link between non-linear Geostatistical methods and non-linear SPDEs. Unfortunately, the state-of-the-art of the theory of non-linear SPDEs seems to be quite technical and sophisticated, and it is not clear how to proceed to create systematically this link with enough generality and simplicity. In this context, this work tries to propose a humble beginning into this aim: constructing an stochastic calculus based on covariance and cross-covariance structures, the basic Geostatistical tools, in a general enough and simple manner. Note that the \textit{basis} of any non-linear stochastic analysis is the definition of \textit{products} between stochastic objects (the stochastic integral being a fundamental case). Those products are the ones involved in the definition of a non-linear SPDE and on the expression of its potential solutions. Hence, exploring the link between stochastic integrals and cross-covariance structures may open the road for a link between Spatial Statistics techniques and the analysis of non-linear SPDEs.

Because all these reasons, we consider that the approach of defining stochastic integrals by simple Riemann sums, with the generality of mean-square continuous stochastic processes as integrands and random measures as integrators, and with their cross-covariance structure determining the possibility of defining such integral, is worth of our work and time.

\subsection{Structure of the paper}
\label{Sec:StructurePaper}

This paper is organized as follows. In Section \ref{Sec:NotationAndPreliminaries} we introduce the basic objects we use in this work which rely on measure theoretical results, as well as the notion of random measure we use: a second-order random measure whose covariance is identified to a measure over the double-dimension space. In Section \ref{Sec:FunctionMeasureKernels} we introduce continuous function-measure kernels $K$ and their related total variation kernel $|K|$. We obtain some analytical results for those concerning iterated integrals ans a Fubini Theorem in Section \ref{Sec:LocalBoundedPartialIntegrals}, and in Section \ref{Sec:CrossPosDefKernels} we study the properties of cross-covariance kernels. Section \ref{Sec:SelfIntegral} is devoted to the definition of the self-integral of a kernel $K$. It is here proven that for some kernels such notion is not uniquely-defined. An auxiliary concept of second-order kernel $K^{(2)}$ is introduced from a cross-tensor product between $K$ and itself, and we define a weaker self-integral for it, called the \textit{quasi-self-integral}.Section \ref{Sec:SelfIntegrabilityAndStochasticIntegration}  is the most important one. It is devoted to the connection between self-integrability and uniquely-defined stochastic integrals. In Section \ref{Sec:FirstApproachStochIntegrals} we define a uniquely-defined stochastic integral between a mean-square continuous process and a random measure over bounded sets. We prove that the expectation of such integral is the self-integral of the cross-covariance kernel between the integrand and the integrator (the self-integrability being hence necessary for unique-definiteness). Section \ref{Sec:GaussianCase} we present our main Theorem \ref{Theo:IntZMGaussian}, which says that when the integrand and the integrator are jointly Gaussian, the unique-definiteness of the stochastic integral is equivalent to the self-integrability of the cross-covariance kernel and the quasi-self-integrability of its second order kernel. Section \ref{Sec:Subsets} focuses on stochastic and self-integrability on subsets (which does not follow immediately from our definitions). It turns out that the self-integral and the stochastic integral are additive functions when uniquely-defined. Theorem \ref{Theo:Resume} synthesize these results in the jointly Gaussian case and presents formulae for the mean and the covariance between stochastic integrals. In Section \ref{Sec:MeasureStructure} the potential measure structure of the self-integral and the stochastic integrals are studied. We give sufficient conditions for those to determine measures over the space. In Section \ref{Sec:Applications} is devoted to some applications of these results. In Section \ref{Sec:TensorProductOfGaussianRandomMeasure} we construct tensor products of jointly Gaussian random measures and we give their mean and covariance structures. Section \ref{Sec:FracBrownianMotion} focuses on the case of fractional Brownian motion $B_{H}$ with Hurst index $H > \frac{1}{2}$, for which its derivative is a random measure \citep{borisov2006constructing} enough well-behaved so the stochastic integral of $B_{H}$ with respect to its derivative is uniquely-defined. In Section \ref{Sec:DerivativeRandomMeasures} we study the case where the derivative of the integrand is an orthogonal random measure without atoms. In such case, the corresponding stochastic integral is never uniquely-defined. This can be interpreted as follows: even if the tensor products random measures can be uniquely-defined, the corresponding double integrals do not have a canonical disintegration. We finalise this paper with the conclusive  Section \ref{Sec:Perspectives}. In Section \ref{Sec:OpenQuestions} we propose some open questions inspired from the developments here presented. In Section \ref{Sec:GeneralizedTensorProducts} we show that the results here obtained can be used to study more multiplicative and tensor products between more general stochastic structures in a manner analogous to the stochastic integral.

\section{Notation and Preliminaries}
\label{Sec:NotationAndPreliminaries}

We denote $\mathbb{N} = \lbrace 0 , 1 , 2 , ... \rbrace$. For a topological space $D$, $C(D)$ denotes the space of complex continuous functions over $D$, and $C_{c}(D)$ the space of complex continuous functions with compact support. If $D \subset \Rd$ is a Borel set, we denote $\Borel{D}$ the $\sigma$-algebra of Borel subsets of $D$, and $\BoundedBorel{D}$ the sub-collection of $\Borel{D}$ consistent of bounded sets (which is a $\delta$-ring in general). We denote $\lambda$ the Lebesgue measure. The space of complex (Borel) measurable functions over $D$ is denoted $\mathcal{M}(D)$.  $\mathcal{M}_{B}(D)$ and $\mathcal{M}_{B,c}(D)$ are the subspaces of $\mathcal{M}(D)$ consisting of bounded functions and bounded compactly-supported functions respectively. $\| \cdot \|_{\infty,A}$ denotes the supremum norm over the set $A$ ($A$ is omitted if it is clear from context). $\mathbbm{1}_{A}$ is the indicator function of the set $A$. LDCT stands for \textit{Lebesgue Dominated Convergence Theorem}.

Concerning probabilistic definitions, we work always with a fixed probability space $\OmegaAP$. A real \textit{stochastic process} $Z$ is a collection of real random variables indexed by a non-empty set $(Z(x))_{x \in D}$ ($D$ is not necessarily a time interval). If $Z(x) \in \LtwoOmega$ for every $x \in D$, then $Z$ is said to be a \textit{second-order} process and  the functions  $m_{Z}(x) = \mathbb{E}(Z(x))$ and $C_{Z}(x,y) = \Cov( Z(x) , Z(y)) $ are called the \textit{mean} and the \textit{covariance functions} respectively. If $D \subset \Rd$, a stochastic process $(Z(x))_{x \in D}$ is said to be mean-square continuous over $D$ if for every $x_{0} \in D$,  $Z(x) \xrightarrow[x \to x_{0}]{L^{2}(\Omega)} Z(x_{0})$. It is known that $Z$ is mean-square continuous if and only if $m_{Z} \in C(D)$ and $C_{Z} \in C(D\times D)$ \citep[Section 21]{sobczyk1991stochastic}.

\subsection{Basic concepts in measure theory}
\label{Sec:BasicConcepts}

We focus on the measurable space $(\Rd , \BorelRd)$ and on subsets of it. 

\begin{Def}
\label{Def:Measure}
A real measure over $D \in \Borel{\Rd}$ is a function $\mu : \BoundedBorel{D} \to \mathbb{R}$ which is $\sigma$-additive over $\BoundedBorel{D}$.
\end{Def} 

We remark that $\mu$ may not be defined over the whole Borel $\sigma$-algebra (some authors would use the word \textit{pre-measure} for $\mu$ \citep{kupka1978caratheodory}). Definition \ref{Def:Measure} implies that $\mu$ is locally finite hence Radon since we are working over the Euclidean space. The real vector space of real measures over $D$ is denoted $\mathscr{M}(D)$. $\mu$ is said to be positive if it takes non-negative values. The \textit{total-variation} measure of $\mu$ is the function $|\mu| : \BoundedBorel{D} \to \R^{+}$ defined by
\begin{equation}
\label{Eq:DefTotalVariation}
|\mu|(A) := \sup \left\lbrace \sum_{j \in J} |\mu(I_{j}) | \ \big| \ (I_{j})_{j \in J} \hbox{ Borel finite partition of } A \right\rbrace.
\end{equation}
$|\mu|$ is indeed a measure \citep[Chapter 6]{rudin1987real} and it is the smallest positive measure such that $|\mu(A)|\leq |\mu|(A)$ for every $A \in \BoundedBorel{D}$. If $D$ is compact the number $|\mu|(D)$ is called the total-variation norm of $\mu$ over $D$, and it defines a norm over $\mathscr{M}(D)$ which turns it into a Banach space.

Let $f \in \mathcal{M}(D)$ and $\mu \in \mathscr{M}(D)$. If $|f|$ is Lebesgue integrable with respect to $|\mu|$ over $D$, then the Lebesgue integral $\int_{D} f d\mu$ is defined, for instance, using the decomposition of $\mu$ into positive and negative parts \citep[Chapter 6]{rudin1987real}. In such case, we use the linear-functional notation for the integral 
\begin{equation}
\label{Eq:BracketNotation}
\langle \mu , f \rangle := \int_{D} f d\mu := \int_{D}f(x)d\mu(x), 
\end{equation}
if $D$ is clear from context. We shall use the bracket or the integral notation depending on the ease of reading of the context. Of course we always have $|\langle \mu , f \rangle | \leq |\mu|(D) \| f \|_{\infty , D}$. 

We can also study measures through the \textit{Bourbakian} approach \citep{bourbaki1965integration}, that is, exploiting the Riesz Representation Theorem interpreting measures as members of dual spaces of continuous functions. We introduce such Theorem in its \textit{topology free} version over the Euclidean space.

\begin{Theo}[\textbf{Riesz Representation}]
\label{Theo:RieszRepresentationRd}
Let $T : C_{c}(\Rd) \to \mathbb{C}$ be a real linear functional. Then, the following are equivalent
\begin{enumerate}[(i)]
\item For every $K \subset \Rd$ compact there exists $C_{K} > 0$ such that
\begin{equation}
\label{Eq:TcontinuousCc}
|\langle T , \varphi \rangle | \leq C_{K} \| \varphi \|_{\infty , K}, \quad \forall \varphi \in C_{c}(\Rd) \hbox{ such that } \supp(\varphi) \subset K.
\end{equation}
\item There exists a unique $\mu \in \mathscr{M}(\Rd)$ such that
\begin{equation}
\label{Eq:TequalsMeasureRiesz}
\langle T , \varphi \rangle = \int_{\Rd} \varphi d\mu, \quad \forall \varphi \in C_{c}(\Rd).
\end{equation}
\end{enumerate}
\end{Theo}

\begin{Rem}
\label{Rem:RieszRepresentationD}
An implication of Theorem \ref{Theo:RieszRepresentationRd} is that if $ D \subset \Rd$ is compact and $C(D)$ is endowed with the supremum norm (hence a Banach space), then $C'(D) = \mathscr{M}(D)$ in a sense analogue to \eqref{Eq:TequalsMeasureRiesz}. The total-variation norm over $\mathscr{M}(D)$ coincides then with the operator norm over $C'(D)$.
\end{Rem}

See \citep[Section 1.10]{tao2022epsilon} for a proof. In this work Theorem \ref{Theo:RieszRepresentationRd} is used mainly to prove that some structures are effectively measures through their identification with linear functionals satisfying \eqref{Eq:TcontinuousCc}.  The following auxiliary definition will be widely used. 

\begin{Def}
\label{Def:RiemannSystem}
Let $ D \in \BoundedBorelRd$. For every $n \in \mathbb{N}$, we consider a finite index set $J_{n}$ such that $|J_{n}|\nearrow \infty$\footnote{$|J|$ denotes the cardinality of the index set $J$.}, a (finite) Borel partition of $D$, $(I_{j}^{n})_{j \in J_{n}}$ satisfying $\displaystyle\max_{j \in J_{n}} \diam( I_{j}^{n}) \xrightarrow[n \to \infty]{} 0 $, and arbitrary tag-points $x_{j}^{n} \in I_{j}^{n}$. The so-constructed sequence $\RiemannSystem{I}{x}{J}$ is called a Riemann system of $D$.
\end{Def}

If $D \in \BoundedBorel{\Rd}$, $f \in C(\Rd)$ and $\mu \in \mathscr{M}(\Rd)$, the Lebesgue integral $\int_{D}fd\mu$ can be computed as the limit of Riemann sums
\begin{equation}
\label{Eq:RiemannSumsintfmu}
\int_{D} f d\mu = \lim_{n \to \infty}\sum_{j \in J_{n}} f(x_{j}^{n})\mu(I_{j}^{n}),
\end{equation}
the limit being independent of the selected Riemann system of $D$ used. This method of constructing integrals is the main focus of this work. Note that any sub-Riemann system, that is, any sub-sequence of the sequence $\RiemannSystem{I}{x}{J}$ is also a Riemann system, and the associated Riemann sums should converge to the same integral.

\subsection{Random measures}
\label{Sec:RandomMeasures}

In the literature it is possible to find many different notions of a \textit{random measure}. Maybe the most intuitive one consists in a random variable taking values in a space of measures. A lot of work has been done around such concept \citep{horowitz1986gaussian,kallenberg2017random,baccelli2020random}. A slightly weaker, but also intuitive and probabilistic-friendly definition is to consider a random variable which takes values in a space of measures almost-surely. This is the case of Poisson and other point processes. However, some important examples of \textit{other kinds} of random measures do not satisfy these requirements. The most important example is Gaussian White Noise over $\Rd$, whose realisations have almost-surely infinite variation over bounded sets \citep[Exercice 3.16]{dalang2009minicourse}.

Another notion of random measure which is often used \citep{morando1969mesures,rao2012random} is the one of a \textit{vector measure} taking values in a (topological) vector space of random variables, such as $L^{p}\OmegaAP$ for $p \in [0,\infty]$. To be precise, if $(E, \mathscr{E} )$ is a measurable space ($\mathscr{E}$ being a $\sigma$-algebra or a $\delta$-ring), a $ L^{p}\OmegaAP$-random measure $M$ over $E$ is a $\sigma$-additive function $M : \mathscr{E} \to L^{p}\OmegaAP$. The most popular case is $p=2$, where $M$ is also called a second-order random measure \citep{thornett1979class,aaberg2011class}. Gaussian White Noise over $\Rd$ enters in such definition. For a second-order random measure, one can define mean and covariance structures given by the functions $A \in \mathscr{E} \mapsto \mathbb{E}(M(A))$ and $(A,B) \in \mathscr{E}\times\mathscr{E} \mapsto  \Cov(M(A),M(B) ) $ respectively.  The first one defines a measure over $(E,\mathscr{E})$, noted $m_{M}$ and called the \textit{mean measure}. The second one defines a  \textit{bi-measure} over $ \mathscr{E}\times\mathscr{E}$, that is, a structure which defines a measure over $A$ when $B$ is fixed and conversely. One important remark is that the covariance structure does not necessarily define a \textit{measure} over the cross-space $( E\times E , \mathscr{E} \otimes \mathscr{E} )$. That is, there are examples of second-order random measures $M$ for which there exists no measure $\Lambda$ over $( E\times E , \mathscr{E} \otimes \mathscr{E} )$ such that $ \Lambda ( A\times B ) = \Cov( M(A) , M(B) )$ \citep[Chapter 2, Example 2]{rao2012random}. Hence some extra precautions must be taken into account when using a second-order random measure. To avoid those, we shall use the following simplified definition. 

\begin{Def}
\label{Def:RandomMeasure}
A second-order measure-covariance random measure (from now on, random measure) over $D \in \Borel{\Rd}$ is a second-order random measure $M : \BoundedBorel{D} \to \LtwoOmega$ such that there exists a measure $C_{M} \in \mathscr{M}(D\times D)$ such that
\begin{equation}
\label{Eq:DefCM}
C_{M}( A \times B ) = \Cov( M(A) , M(B) ), \quad \forall A , B \in \BoundedBorel{D}.
\end{equation}
\end{Def}

Definition \ref{Def:RandomMeasure} is the only sense in which the words ``random measure'' will be used in this work. Some authors have already remarked the convenience of this measure-covariance supposition for different aims \citep{borisov2006constructing,carrizov2022karhunen}. The measure $C_{M}$ is called the \textit{covariance measure} of $M$. The cautious reader will realize that, all along this work, the supposition of $C_{M}$ being a measure will be of great help when proving most part of the results.

If $\varphi \in \mathcal{M}(\Rd)$ satisfies
\begin{equation}
\label{Eq:ConditionsVarphiIntegral}
\langle |m_{M}| , |\varphi |\rangle < \infty, \quad \hbox{and} \quad   \langle |C_{M}| , |\varphi| \otimes |\varphi| \rangle < \infty, 
\end{equation}
then the stochastic integral
\begin{equation}
\label{Eq:IntMvarphi}
\langle M , \varphi \rangle := \int_{\Rd} \varphi dM
\end{equation}
can be defined. This is nothing but the \textit{Dunford-Schwartz integral} \citep[Section IV.10]{dunford1958linear} with respect to the $\LtwoOmega$-valued measure $M$. If $\varphi,\phi$ satisfy \eqref{Eq:ConditionsVarphiIntegral} then  (\citep{borisov2006constructing} or \citep[Proposition 3.3.1]{carrizov2018development})
\begin{equation}
\label{Eq:MeanCovIntegralM}
\mathbb{E}( \langle M , \varphi \rangle ) = \langle m_{M} , \varphi \rangle \quad \hbox{and} \quad \Cov\left( \langle M , \varphi \rangle \langle M , \phi \rangle  \right) = \langle C_{M} , \varphi \otimes \overline{\phi} \rangle.
\end{equation}

Similarly as the case of deterministic measures, random measures can be studied through the \textit{Bourbakian approach}. The following Theorem is a stochastic version of Riesz Representation Theorem \ref{Theo:RieszRepresentationRd} and it is a fundamental tool for such developments. It states that a random measure can be defined when studying only its action against continuous functions, and the integrals with respect to more general measurable functions can be hence defined by extension arguments. We will state in the setting of $C(D)$ with $D$ compact\footnote{It can be extended to the case $C_{c}(D)$ following simple arguments which we omit.}.

\begin{Theo}
\label{Theo:ExtensionRandomMeasure}
Let $D \subset \Rd$ compact. Let $M : C(D) \to \LtwoOmega$ be such that there exist measures $m_{M} \in \mathscr{M}(D)$ and  $C_{M} \in \mathscr{M}(D\times D)$ such that \eqref{Eq:MeanCovIntegralM} hold for every $\varphi,\phi \in C(D)$. Then, there exists a unique random measure $\tilde{M} : \Borel{D} \to  \LtwoOmega$  such that
\begin{equation}
\label{Eq:M=tildeMinC(D)}
\langle M , \varphi \rangle = \int_{D}\varphi d\tilde{M}, \quad \forall \varphi \in C(D).
\end{equation}
\end{Theo}

Of course, the mean and covariance measures of $\tilde{M}$ in Theorem \ref{Theo:ExtensionRandomMeasure} are $m_{M}$ and $C_{M}$, respectively. This Theorem also implies that the definition of $M$ can be extended to $\mathcal{M}_{B}(D)$ linearly and continuously (with the supremum norm). While this Theorem is powerful, it is not the main focus of this work, and it will be used in an auxiliary manner. Its proof is hence refereed  to Appendix \ref{Sec:AppendixProofTheoRieszRandom}.

We make explicit two examples of random measures in the sense of Definition \ref{Def:RandomMeasure}.

\begin{Exa}
\label{Ex:GaussianWhiteNoise}
Gaussian White Noise $W$ over $\Rd$ is a centred Gaussian stochastic process $(W(A))_{A \in \BoundedBorel{\Rd}}$ such that $\Cov(W(A),W(B)) = \lambda(A\cap B)$. Note that $C_{W}$ is the measure related to the linear functional $\psi \in C_{c}(\RdxRd) \mapsto \langle C_{W} , \psi \rangle := \int_{\Rd}\psi(x,x)dx$. The stochastic integral $\langle W , \varphi \rangle$ can actually be extended to $\varphi \in L^{2}(\Rd)$ since in such case $\langle C_{W} , \varphi\otimes \overline{\varphi} \rangle = \int_{\Rd}|\varphi|^{2}dx < \infty$ and thus condition \eqref{Eq:MeanCovIntegralM} hold. $\square$
\end{Exa}

\begin{Exa}
\label{Ex:IntZMuRandomMeasure}
Consider a mean-square continuous process $Z$ over $\Rd$, and let $\mu \in \mathscr{M}(\Rd)$. Then, one can define without ambiguity through Riemann sums (\citep[Theorem 4.5.2]{soong1973random} or \citep[3.2.1]{carrizov2018development}) the application $A \in \BoundedBorel{\Rd} \mapsto \int_{A}Z(x)d\mu(x)$. One has  
\begin{equation}
\label{Eq:MeanCovarianceIntZMu}
\mathbb{E}\left( \int_{A}Z(x)d\mu(x) \right) = \int_{A}m_{Z}d\mu \quad ; \quad \Cov\left(  \int_{A}Z(x)d\mu(x) , \int_{B}Z(x)d\mu(x) \right) = \int_{A\times B} C_{Z} d(\mu\otimes \mu).
\end{equation}
Thus, $A \mapsto \int_{A}Zd\mu $ has mean and covariance structures determined by measures, hence it is a random measure over $\Rd$. $\square$
\end{Exa}

Other examples which will be treated later are fractional Brownian motion with Hurst index $H > \frac{1}{2}$ (Section \ref{Sec:FracBrownianMotion}) and orthogonal random measures (Section \ref{Sec:DerivativeRandomMeasures}).

\section{Function-measure kernels}
\label{Sec:FunctionMeasureKernels}

In this Section we study basic analytical properties and operations on function-measure kernels and we give the first insights on their importance on stochastic integration. We begin with the main definitions.

\begin{Def}
\label{Def:Kernel}
Let $K : \Rd \times \BoundedBorelRp \to \R$ be a mapping such that $K(x,\cdot) \in \mathscr{M}(\R^{p})$ for all $x \in \Rd$, and $K(\cdot , A) \in C(\Rd)$ for all $A \in \BoundedBorelRp$.
Then, $K$ is called a continuous-function$-$measure kernel (from now on, cf-m kernel) over $\Rd\times\R^{p}$.  
\end{Def}

In the rest of this Section $K$ will always denote a cf-m kernel over $\Rd\times\R^{p}$. 

\begin{Def}
\label{Def:TotalVariationKernel}
The total-variation kernel of $K$ is the kernel $|K| : \Rd \times \BoundedBorel{\R^{p}}\to \R^{+}$ given by $|K|(x,A) := |K(x , \cdot)|(A)$.
\end{Def}

$|K|$ can be computed using formula \eqref{Eq:DefTotalVariation}. Note that if $(I_{j})_{j}$ is a finite partition of $A$, the sum $\sum_{j} |K(\cdot , I_{j})| $ defines a continuous function. $|K|(\cdot,A)$ is hence the supremum of a set of continuous functions, so it is lower semi-continuous and thus measurable. However, it is not clear if it is locally bounded, which is a regularity property we shall ask and study in the following Section.

\subsection{Local boundedness and partial integrals}
\label{Sec:LocalBoundedPartialIntegrals}

\begin{Def}
\label{Def:KlocallyBounded}
$K$ is said to be locally bounded if
\begin{equation}
\label{Eq:DefKLocallyBounded}
\sup_{x \in A} |K|(x,B) < \infty, \quad \forall A \in \BoundedBorelRd, B \in \BoundedBorelRp.
\end{equation}
In such case we define $\| K \|_{A,B,\varepsilon} := \sup_{x \in A} |K|(x,B)$.
\end{Def}

The notation $\| K \|_{A,B,\varepsilon}$ is inspired from the injective norm of a tensor product of normed spaces, but we shall not give much attention to this. The extra requirement of local boundedness allows to gain some regularity when operating with $K$. Note that integrals of the form $\int_{\R^{p}}\psi(x,y)dK(x,\cdot)(y)$ are well-defined for every $\psi \in \mathcal{M}_{B,c}(\Rd\times\R^{p})$ and every $x$.

\begin{Prop}
\label{Prop:IntPsiKyContinuous}
Let $\psi \in C_{c}(\Rd\times\mathbb{R}^{p})$. If $K$ is locally bounded, then the function $x \mapsto \int_{\mathbb{R}^{p}} \psi(x,y) dK(x,\cdot)(y)$ is continuous with compact support.
\end{Prop}

\textbf{Proof:} Suppose $\supp(\psi) \subset D_{d}\times D_{p}$ with $D_{d} \subset \Rd$ and $D_{p} \subset \mathbb{R}^{p}$ compacts. If $\RiemannSystem{I}{y}{J}$ is a Riemann system of $D_{p}$, then we have the point-wise limit of continuous functions
\begin{equation}
\label{Eq:IntPsiKyRiemannSum}
\int_{\mathbb{R}^{p}} \psi(x,y) dK(x,\cdot)(y) = \lim_{n \to \infty} \sum_{j \in J_{n}} \psi( x , y_{j}^{n}) K(x , I_{j}^{n}), \quad \forall x \in \Rd.
\end{equation}
\eqref{Eq:IntPsiKyRiemannSum} is null for $x \notin D_{d}$, hence the limit has compact support. Let us prove the convergence is uniform. Let $\epsilon > 0 $. Since $\psi$ is uniformly continuous over $D_{d}\times D_{p}$, there exists $\delta > 0 $ such that if $|x-u|+|y-v| \leq \delta$, $|\psi(x,y)-\psi(u,v)| \leq \epsilon \| K \|_{D_{d},D_{p},\varepsilon}^{-1}$. If $n$ is big enough so $\max_{j \in J_{n}}\diam( I_{j}^{n}) \leq \delta $, then
\small
\begin{equation}
\label{Eq:ProofUnifConvergenceIntPsiKy}
\begin{aligned}
\sup_{x \in D_{d}}\left| \int_{D_{p}} \psi( x , y ) d K(x,\cdot)(y) - \sum_{j \in J_{n}} \psi( x , y_{j}^{n})K(x, I_{j}^{n}) \right| &= \sup_{x \in D_{d}} \left| \sum_{j \in J_{n}} \int_{I_{j}^{n}} \left( \psi(x,y)-\psi(x,y_{j}^{n}) \right)d K(x , \cdot) (y)   \right| \\
&\leq \sup_{x \in D_{d}} \sum_{j \in J_{n}} \underbrace{\sup_{y \in I_{j}^{n}} |\psi( x, y) - \psi( x, y_{j}^{n} ) |}_{\leq \epsilon \|K \|_{D_{d},D_{p},\varepsilon}^{-1} }  |K|(x,I_{j}^{n}) \\
&\leq \sup_{x \in D_{d}} \epsilon \|K \|_{D_{d},D_{p},\varepsilon}^{-1} \underbrace{\sum_{j \in J_{n}} |K|(x,I_{j}^{n})}_{= |K|(x ,D_{p})}   = \epsilon
\end{aligned} 
\end{equation}
\normalsize
This proves a uniform convergence hence the limit function \eqref{Eq:IntPsiKyRiemannSum} is continuous. $\blacksquare$

A similar regularity result holds when $\psi$ is independent of $x$.
\begin{Prop}
\label{Prop:IntPhiKyContinuous}
Let $\phi \in \mathcal{M}_{B,c}(\R^{p})$. If $K$ is locally bounded, then the function $x \mapsto \int_{\R^{p}} \phi(y) dK(x,\cdot)(y)$ is continuous.
\end{Prop}

\textbf{Proof:} Let $D_{p} \subset \R^{p}$ compact with $\supp(\phi) \subset D_{p}$. Since $\phi$ is bounded it can be expressed as a uniform-over-$D_{p}$ limit of simple functions $(\phi_{n})_{n}$. For a simple function $\phi_{n}$ it is clear that $x \mapsto  \int_{\R^{p}} \phi_{n}(y) dK(x,\cdot)(y)$ is continuous. Now, one has for every compact $D \subset \Rd$,
\begin{equation}
\label{Eq:ProofIntPhiKyContinuous}
\sup_{x \in D}\left|\int_{D_{p}} \phi(y) - \phi_{n}(y) dK(x , \cdot)(y) \right| \leq \| K \|_{D,D_{p},\varepsilon} \| \phi - \phi_{n} \|_{\infty} \xrightarrow[n \to \infty]{} 0.
\end{equation}
The function $x \mapsto \int_{\R^{p}} \phi(y) dK(x,\cdot)(y)$ is hence continuous as a uniform-on-compacts limit of continuous functions. $\blacksquare$

Consider now for any given $\mu \in \mathscr{M}(\Rd)$ the iterated integral
\begin{equation}
\label{Eq:IntPsiKxy(dp)}
\int_{\Rd}\int_{\R^{p}} \psi(x,y) dK(x,\cdot)(y) d\mu (x),
\end{equation}
which is well-defined for $\psi \in C_{c}(\Rd\times\R^{p})$. This induces a  real linear functional from $C_{c}(\Rd\times\mathbb{R}^{p})$ to $\mathbb{C}$.
\begin{Prop}
\label{Prop:ExistenceMeasureKmu}
Let $\mu \in \mathscr{M}(\Rd)$. If $K$ is locally bounded, there exists a unique $\tilde{K}_{\mu} \in \mathscr{M}(\Rd\times\mathbb{R}^{p})$ such that
\begin{equation}
\label{Eq:KmuMeasureEqualsK}
\int_{\Rd}\int_{\R^{p}} \psi(x,y) dK(x,\cdot)(y) d\mu (x) = \langle \tilde{K}_{\mu} , \psi \rangle, \quad \forall \psi \in C_{c}(\Rd\times\mathbb{R}^{p}).
\end{equation} 
\end{Prop}

\textbf{Proof:} If $\supp(\psi) \subset D_{d}\times D_{p} \subset \Rd\times\mathbb{R}^{p}$, $D_{d}$ and $D_{p}$ compact, then
\begin{equation}
\label{Eq:ProofKDefinesContLinearFunctional}
\left| \int_{\Rd}\int_{\R^{p}} \psi(x,y) dK(x,\cdot)(y) d\mu (x)  \right| \leq |\mu|(D_{d}) \| K \|_{D_{d},D_{p},\varepsilon} \| \psi \|_{\infty}.
\end{equation}
From Riesz Representation Theorem \ref{Theo:RieszRepresentationRd}, there exists a unique measure $\tilde{K}_{\mu}$ identifying the linear functional $\psi \mapsto \int_{\Rd}\int_{\R^{p}} \psi(x,y) dK(x,\cdot)(y) d\mu (x)$ . $\blacksquare$

It follows that the iterated integral \eqref{Eq:IntPsiKxy(dp)} can be extended uniquely to any function in $\mathcal{M}_{B,c}(\Rd\times\R^{p})$ through the use of the measure $\tilde{K}_{\mu}$ having hence 
\begin{equation}
\label{Eq:K=Kmu(dp)}
\int_{\Rd} \int_{\R^{p}} \psi( x,y ) dK(x,\cdot)(y) d\mu(x) = \int_{\Rd\times \R^{p}} \psi(x,y) d\tilde{K}_{\mu} (x,y), \quad \forall \psi \in \mathcal{M}_{B,c}(\Rd\times \R^{p} ). 
\end{equation}
Note that for a fixed $\psi \in \mathcal{M}_{B,c}(\Rd\times \R^{p} )$, the function $x \mapsto \int_{\R^{p}}\psi(x,y)dK(x,\cdot)(y)$, which is compactly supported and bounded by $\| \psi \|_{\infty}\| K \|_{D_{d},D_{p},\varepsilon}$, is nothing but the Radon-Nikodym derivative of the measure $A \in \BoundedBorel{\Rd} \mapsto \int_{A\times \R^{p}} \psi(x,y) d\tilde{K}_{\mu}(x,y)$ with respect to $\mu$. It is hence in $\mathcal{M}_{B,c}(\Rd)$.

Let us now study the iterated integral ``\textit{in the other sense}'', first integrating  over $\Rd$ and then over $\R^{p}$. This would mean to give sense to the integral
\begin{equation}
\label{Eq:IntPsiKxy(pd)(informal)}
``\int_{\R^{p}} d\left( \int_{\Rd} \psi(x,\cdot) K(x,\cdot) d\mu(x) \right)(y)".
\end{equation}
This notation makes no further sense until now for a general $\psi$. Note that the second component of $\psi$ acts as part of the measure integrating over $\R^{p}$. One option for giving a meaning to this is to use a disintegration of the measure $\tilde{K}_{\mu}$ in the opposite order that in \eqref{Eq:K=Kmu(dp)}. Rather than doing this, we will consider a direct construction of the iterated integrals in the new order, starting from simple cases of $\psi$ and then extending. Let $\varphi \in \mathcal{M}_{B,c}(\Rd)$ and $\mu \in \mathscr{M}(\Rd)$. We consider the object
\begin{equation}
\label{Eq:IntVarphiKmux(informal)}
``\int_{\Rd}\varphi(x)K(x,\cdot)d\mu(x)".
\end{equation}
Intuitively, this integral should define a measure over $\R^{p}$. Indeed, the application $ A \in \BoundedBorelRp \mapsto \int_{\Rd}\varphi(x)K(x,A)d\mu(x) $ is clearly well-defined and additive.

\begin{Prop}
\label{Prop:IntVarphiKmuxDefinesMeasure}
If $K$ is locally bounded, $ A \mapsto \int_{\Rd}\varphi(x)K(x,A)d\mu(x) $ defines a measure over $\R^{p}$.
\end{Prop}

\textbf{Proof:} We need to prove the $\sigma$-additivity. Consider a sequence of bounded Borel sets of $\R^{p}$, $(D_{n})_{n \in \mathbb{N}}$, with $D_{n} \searrow \emptyset$. Set the functions $f_{n} := K ( \cdot ,D_{n})$. $(f_{n})_{n \in \mathbb{N}}$ converges point-wise to $0$ and dominated by $ x \mapsto |K|(x,D_{0})$, which is locally bounded since $K$ is locally bounded. From LDCT it follows
\begin{equation}
\label{Eq:LDCTinProofIntVarphiKmuxMeasure}
\lim_{n \to \infty} \int_{\Rd} K(x,D_{n})\varphi(x) d\mu(x) = \int_{\supp(\varphi)} \lim_{n \to \infty} f_{n}(x)\varphi(x) d\mu(x) = 0.
\end{equation}
This proves the $\sigma-$additivity of $ A \mapsto \int_{\Rd}\varphi(x)K(x,A)d\mu(x) $. $\blacksquare$

Since \eqref{Eq:IntVarphiKmux(informal)} is a measure, the iterated integral
\begin{equation}
\int_{\R^{p}}\phi(y) d\left( \int_{\Rd}K(x,\cdot)\varphi(x)d\mu(x) \right)(y), \quad \varphi \in \mathcal{M}_{B,c}(\Rd), \phi \in \mathcal{M}_{B,c}(\R^{p})
\end{equation}
makes sense. By linearity, we can give a meaning to \eqref{Eq:IntPsiKxy(pd)(informal)} for every $\psi$ in the tensor product space $\mathcal{M}_{B,c}(\Rd)\otimes\mathcal{M}_{B,c}(\R^{p}) $. Now, we ask ourselves if such iterated integral coincides with the iterated integral \eqref{Eq:IntPsiKxy(dp)}. We begin with the continuous case. 

\begin{Prop} 
\label{Prop:FubiniKPhiVarphiContinuous}
Let $\mu \in \mathscr{M}(\Rd)$, $\varphi \in C_{c}(\Rd)$ and $\phi \in C_{c}(\R^{p})$. If $K$ is locally bounded, then
\begin{equation}
\label{Eq:FubiniKPhiVarphiContinuous}
\int_{\Rd}\varphi(x)\int_{\R^{p}} \phi(y) dK(x,\cdot)(y) d\mu(x)  = \int_{\R^{p}}\phi(y) d\left( \int_{\Rd}K(x,\cdot)\varphi(x)d\mu(x) \right)(y).
\end{equation}
\end{Prop}

\textbf{Proof:}  Let $D_{d},D_{p}$ compact with $\supp(\varphi) \subset D_{d}$ and $\supp(\phi) \subset D_{p}$. Let $ \RiemannSystem{I}{y}{J}$ be a Riemann system of $D_{p}$. Since the convergence $\sum_{j \in J_{n}} \phi( y_{j}^{n}) K(x , I_{j}^{n}) \to \int_{K_{p}} \phi(y) dK(x , \cdot)(y) $ is uniform over $x \in D_{d}$ (proof of Proposition \ref{Prop:IntPsiKyContinuous}), we can interchange integral and limit symbols:
\begin{equation}
\begin{aligned}
\int_{D_{p}} \phi(y) d\left( \int_{D_{d}} \varphi(x) K( x , \cdot ) d\mu(x) \right) (y) &= \lim_{n \to \infty} \sum_{j \in J_{n}} \phi( y_{j}^{n}) \int_{D_{d}} \varphi(x) K( x , I_{j}^{n}) d\mu(x) \\
&=  \int_{D_{d}} \varphi(x)\lim_{n \to \infty}  \sum_{j \in J_{n}} \phi( y_{j}^{n} ) K( x , I_{j}^{n} )   d\mu(x) \\
&= \int_{D_{d}} \varphi(x) \int_{D_{p}} \phi(y) dK(x, \cdot )(y) d\mu(x). \quad \blacksquare
\end{aligned}
\end{equation}

By linearity, for all $\psi \in C_{c}(\Rd)\otimes C_{c}(\R^{p})$ we have
\begin{equation}
\label{Eq:FubiniPsiKmu}
\int_{\R^{p}} d\left( \int_{\Rd}\psi(x,\cdot)K(x,\cdot)d\mu(x) \right)(y) =  \int_{\Rd}\int_{\R^{p}} \psi(x,y) dK(x,\cdot)(y) d\mu(x)  =  \langle \tilde{K}_{\mu} , \psi \rangle.
\end{equation}
But since $\tilde{K}_{\mu}$ is a measure, its definition can be extended uniquely to any $\psi \in C_{c}(\Rd \times \R^{p})$ using the density of $ C_{c}(\Rd)\otimes C_{c}(\R^{p}) $ in $C_{c}(\Rd \times \R^{p})$ with the uniform-on-compacts topology. The extension of the three integrals in \eqref{Eq:FubiniPsiKmu} to any $\psi \in \mathcal{M}_{B,c}(\Rd\times\R^{p})$ has thus a unique meaning and they all must agree, giving a complete rigorous sense to the integral \eqref{Eq:IntPsiKxy(pd)(informal)}. We conclude the following Theorem.

\begin{Theo}[\textbf{Fubini}]
\label{Theo:FubiniKernelMu}
If $K$ is locally bounded, then for every $\psi \in \mathcal{M}_{B,c}(\Rd\times\R^{p})$ the iterated integral \eqref{Eq:IntPsiKxy(pd)(informal)} is well-defined and the equalities \eqref{Eq:FubiniPsiKmu} hold.
\end{Theo}

As we shall see later in Section \ref{Sec:SelfIntegral}, the integral we have defined in \eqref{Eq:IntPsiKxy(pd)(informal)} is the first non-trivial example of \textit{self-integral} we have used in this work.

\subsection{Cross-positive-definite kernels}
\label{Sec:CrossPosDefKernels}

The most important kind of cf-m kernels over $\Rd\times \R^{p}$ we use in this work are those of the type
\begin{equation}
\label{Eq:DefKZM}
K_{Z,M}(x,A) := \Cov(Z(x) , M(A)),
\end{equation}
$Z$ being a mean-square continuous stochastic process over $\Rd$ and $M$ being random measure over $\R^{p}$. The fact that $K_{Z,M}$ is indeed a cf-m kernel can be deduced from the continuity of $Z$ and the $\sigma$-additivity of $M$. The kernel $K_{Z,M}$ is called the \textit{cross-covariance kernel} between $Z$ and $M$. We shall see that such kind of kernel has particular extra regularity properties.

We recall that the covariance function $C_{Z}$ and the covariance measure $C_{M}$ must have a positive-definite behaviour. To be precise, $C_{Z}$  is such that if $(x_{1} , ... , x_{n}) \in (\Rd)^{n} $ and if $(\alpha_{1} , ... , \alpha_{n}) \in \mathbb{C}^{n}$, then
$\sum_{j,k=1}^{n} \alpha_{j} C_{Z}(x_{j} , x_{k} ) \overline{ \alpha_{k}} \geq 0.$ Similarly, $C_{M}$ is such that if $(E_{1} , ... , E_{n}) \in (\mathcal{B}_{B}(\R^{p}) )^{n} $ and if $(\alpha_{1} , ... , \alpha_{n}) \in \mathbb{C}^{n}$, then $\sum_{j,k=1}^{n} \alpha_{j} C_{M}(E_{j} \times E_{k} ) \overline{ \alpha_{k}} \geq 0$. $C_{Z}$ and $C_{M}$ describe roughly the statistical properties of $Z$ and $M$ separately, but $K_{Z,M}$ describes, also roughly in general, their dependence interaction. Based on this, we have the following definition and criterion of cross-positivity.

\begin{Def}
\label{Def:CrossPosDefK}
We say that $K$ is a \textit{cross-positive-definite kernel} (abbreviated cross-pos-def kernel), or that it is a \textit{cross-covariance kernel}, if there exist a covariance function $C_{Z} \in C(\RdxRd)$ and a covariance measure $C_{M} \in \mathscr{M}(\R^{p} \times \R^{p})$ such that 
\begin{equation}
\label{Eq:DefCrossPosDefK}
\left|\sum_{j=1}^{n} \sum_{k=1}^{m} \alpha_{j} \beta_{k} K(x_{j} , E_{k} ) \right| \leq \frac{1}{2}\left[ \sum_{i,j=1}^{n}\alpha_{j}C_{Z}(x_{i},x_{j})\overline{\alpha_{k}}   +   \sum_{k,l=1}^{m}\beta_{k}C_{M}( E_{k} \times E_{l} )\overline{\beta_{l}} \right].
\end{equation}
for any $n,m \geq 1$, $(\alpha_{1} , ... , \alpha_{n}, \beta_{1} , ... , \beta_{m}) \in \mathbb{C}^{n+m}$, $(x_{1} , ... , x_{n} ) \in (\Rd)^{n}$ and $(E_{1} , ... , E_{m} ) \in \mathcal{B}_{B}(\R^{p})^{m}$. 
\end{Def}

\begin{Prop}
\label{Prop:CrossPosDefIsCrossCov}
$K$ is cross-pos-def if and only if $K(x,A) = \Cov(Z(x) , M(A) )$ for a mean-square continuous process $Z$ over $\Rd$ and a random measure $M$ over $\R^{p}$.
\end{Prop}

\textbf{Proof:} For the sufficiency, condition \eqref{Eq:DefCrossPosDefK} is actually equivalent to the non-negative-definiteness of the covariance matrix of the random vector $( Z(x_{1}) , ... , Z(x_{n}) , M(E_{1}) , ... , M(E_{m}) )$. For the necessity, one can use the classical Kolmogorov existence Theorem to construct a probability space where $Z$ and $M$ are well-defined jointly Gaussian processes, with covariances $C_{Z}$ and $C_{M}$ respectively, and cross-covariance $K$. Condition \eqref{Eq:DefCrossPosDefK} guarantees the Kolmogorov compatibility conditions when constructing the covariance matrices of finite-dimensional Gaussian random vectors. Details are left to the reader. $\blacksquare$

The first self-evident (we omit the proof) consequence to is the following.

\begin{Prop}[\textbf{Cauchy-Schwarz inequality}]
\label{Prop:CauchySchwarz}
If $K$ is cross-pos-def, and $C_{Z}, C_{M}$ are as in Definition \ref{Def:CrossPosDefK}, then
\begin{equation}
\label{Eq:CauchySchwarzK}
\left|   \sum_{j=1}^{n} \sum_{k=1}^{m} \alpha_{j} \beta_{k} K(x_{j} , E_{k} )  \right| \leq \sqrt{  \sum_{i,j=1}^{n}\alpha_{j}C_{Z}(x_{i},x_{j})\overline{\alpha_{k}}  } \sqrt{ \sum_{k,l=1}^{m}\beta_{k}C_{M}( E_{k} \times E_{l} )\overline{\beta_{l}}}, 
\end{equation}
for any $n,m \geq 1$, $(\alpha_{1} , ... , \alpha_{n}, \beta_{1} , ... , \beta_{m}) \in \mathbb{C}^{n+m}$, $(x_{1} , ... , x_{n} ) \in (\Rd)^{n}$ and $(E_{1} , ... , E_{m} ) \in \mathcal{B}_{B}(\R^{p})^{m}$.
\end{Prop}

The first really important result on cross-pos-def kernels is the following.

\begin{Theo}
\label{Theo:CrossPosDefKisLocBound}
If $K$ is cross-pos-def, then it is locally bounded.
\end{Theo}

\textbf{Proof: } Let $A \in \mathcal{B}_{B}(\R^{p})$ and $D \in \BoundedBorelRd$. Then, by definition
\begin{equation}
\sup_{x \in D} |K|(x , A) = \sup_{x \in D} \sup_{(E_{j})_{j} \in \mathscr{P}_{F}(A) }  \sum_{j}|K(x,E_{j})|, 
\end{equation}
where $\mathscr{P}_{F}(A)$ denotes the collection of all finite Borel partitions of $A$. For any $x \in D$ and any $(E_{j})_{j} \in \mathscr{P}_{F}(A)$, we set $\alpha_{x,j} := \sgn( K(x,E_{j} ) )$. These coefficients satisfy hence $|\alpha_{x,j}| = 1$ and $ |K(x,E_{j})| = \alpha_{x,j}K(x,E_{j}) $. From Cauchy-Schwarz inequality \eqref{Eq:CauchySchwarzK}, we obtain
\begin{equation}
\label{Eq:|K|LocallyBoundedProof}
\begin{aligned}
\sup_{x \in D} |K|(x ,A) &= \sup_{x \in D} \sup_{(E_{j})_{j} \in \mathscr{P}_{F}(A) }\left|  \sum_{j}\alpha_{j,x}K(x,E_{j}) \right| \\
&\leq \sup_{x \in D} \sup_{(E_{j})_{j} \in \mathscr{P}_{F}(A) } \sqrt{C_{Z}(x,x)} \sqrt{\sum_{j,k} \alpha_{j,x} C_{M}(E_{j}\times E_{k}) \alpha_{k,x} } \\
& \leq \sup_{x \in D} \sqrt{C_{Z}(x,x)}\sup_{(E_{j})_{j} \in \mathscr{P}_{F}(A) } \sqrt{\sum_{j,k} |C_{M}|(E_{j}\times E_{k})  } \\
& \leq \sup_{x \in D} \sqrt{C_{Z}(x,x)} \sqrt{|C_{M} |(A\times A)} < \infty,
\end{aligned} 
\end{equation} 
where we have used that $C_{Z}$ is continuous and $C_{M}$ is a measure. $\blacksquare$

Results obtained in Section \ref{Sec:LocalBoundedPartialIntegrals} apply thus for cross-pos-def kernels. New surprises are also present in this case. As said above, it is not clear in general if the function $|K|(\cdot , A)$ is continuous. It turns out that for cross-pos-def kernels an even stronger regularity property holds. To present it, we will introduce the following auxiliary definition.

\begin{Def}
\label{Def:IncrementKernel}
We define the increment kernel of $K$ as the kernel $K_{I}$ over $(\RdxRd)\times \R^{p}$ given by $K_{I}( (x,y) , A ) := K(x,A) - K(y,A)$.
\end{Def}

If $K$ is cross-pos-def, then $K_{I}$ is also cross-pos-def. This because if $K$ is the cross-covariance kernel between $Z$ and $M$, then $K_{I}$ is the cross-covariance kernel between the increment process
\begin{equation}
\label{Eq:IncrementProcess}
I_{Z}(x,y) := Z(x) - Z(y)
\end{equation}
and $M$. We will use the notation $K_{I_{Z} , M}$ for the increment kernel of $K_{Z,M}$ in the cross-covariance case. Note that the covariance function of $I_{Z}$ is given by
\begin{equation}
\label{Eq:CovIncrementProcess}
C_{I_{Z}}( (x,y) , (u,v) ) := C_{Z}(x,u)- C_{Z}(x,v) - C_{Z}(y,u) + C_{Z}(y,v),
\end{equation}
which is continuous over $\Rd\times\Rd\times\Rd\times\Rd$, hence $I_{Z}$ is mean-square continuous. 

\begin{Theo}
\label{Theo:KcrossPosDefEquicontinuous}
Let $D \subset \Rd$ compact and $A \in \BoundedBorel{\R^{p}}$. If $K$ is cross-pos-def, then for every $\epsilon > 0 $ there exists $\delta > 0 $ such that 
\begin{equation}
\label{Eq:IncrementsEquiBounded}
\sup_{\stackrel{x,y \in D}{|x-y|\leq \delta}}\sup_{B \in \mathcal{B}(A)} | K(x , \cdot ) - K(y , \cdot)|(B) \leq \epsilon.
\end{equation}
In particular, the family of functions $\lbrace |K|(\cdot,B) \ \big| \ B \in \mathcal{B}(A)  \rbrace$ is uniformly equicontinuous over $D$.
\end{Theo}

\textbf{Proof: } Let $\epsilon > 0 $. We note $K = K_{Z,M}$. Since the increment kernel $K_{I_{Z},M}$ is cross-pos-def, we can use Cauchy-Schwarz inequality \eqref{Eq:CauchySchwarzK}. If $B \in \mathcal{B}(A)$, then
\begin{equation}
\label{Eq:CauchySchwarzIncrements}
\begin{aligned}
\sup_{ x,y \in D} |K_{I_{Z},M}((x,y) , \cdot ) |(B) &\leq \sup_{ x,y \in D } \sqrt{C_{I_{Z}}( (x,y),(x,y) ) } \sqrt{ |C_{M}|(B \times B )} \\
&\leq \sup_{ x,y \in D } \sqrt{C_{I_{Z}}( (x,y),(x,y) ) } \sqrt{ |C_{M}|(A \times A )}. 
\end{aligned}
\end{equation} 
Since $C_{I_{Z}}$ is continuous and $D$ compact, we can use the uniform continuity of $C_{I_{Z}}$ over $D\times D\times D\times D$ and take $\delta > 0 $ such that if $|x-y|\leq\delta$, then $C_{I_{Z}}((x,y),(x,y)) \leq \epsilon^{2} \frac{1}{|C_{M}|(A\times A)}$. The inequality \eqref{Eq:IncrementsEquiBounded} comes then immediately from \eqref{Eq:CauchySchwarzIncrements}. The uniform equicontinuity of the family $\lbrace |K|(\cdot,B) \ \big| \ B \in \mathcal{B}(A)  \rbrace$ comes from $| |K_{Z,M}|(x,B) - |K_{Z,M}|(y , B) | \leq |K_{Z,M}(x,\cdot) - K_{Z,M}(y,\cdot) |(B)$, which is just the triangular inequality applied to the total-variation norm over $B$. $\blacksquare$

We finally consider the iterated integrals studied in Section \ref{Sec:LocalBoundedPartialIntegrals} and their relation to $Z$ and $M$ in the case of a cross-pos-def kernel.

\begin{Prop}
\label{Prop:IteratedIntegralKcrossPos}
Let $K_{Z,M}$ be the cross-covariance kernel between $Z$ and $M$. Then, for every $\mu \in \mathscr{M}(\Rd), \varphi \in \mathcal{M}_{B,c}(\Rd)$, and $ \phi \in \mathcal{M}_{B,c}(\R^{p})$, one has
\begin{equation}
\label{Eq:IteratedIntegralKCrossPos}
\int_{\Rd} \int_{\R^{p}}\varphi(x)\overline{\phi}(y) dK_{Z,M}(x,\cdot)(y) d\mu(x) = \mathbb{C}ov\left( \int_{\Rd} Z(x) \varphi(x) d\mu(x) , \int_{\R^{p}} \phi(y) dM(y) \right).
\end{equation}
\end{Prop}

\textbf{Proof:} Consider first the case $\phi = \mathbbm{1}_{A}$, $A \in \BoundedBorel{\R^{p}}$. Let $\RiemannSystem{I}{x}{J}$ be a Riemann system of $\supp(\varphi)$. Then,
\begin{equation}
\label{Eq:ProofIteratedIntKCrossPosM(A)}
\begin{aligned}
\Cov\left(  \int_{\Rd} \varphi(x)Z(x)d\mu(x) ,  M(A)  \right) &= \lim_{n \to \infty} \Cov\left( \sum_{j \in J_{n}}   Z(x_{j}^{n}) \int_{I_{j}^{n}} \varphi d\mu \ , \  M(A)  \right) \\
&= \lim_{n \to \infty}  \sum_{j \in J_{n}} K_{Z,M}( x_{j}^{n} , A ) \int_{I_{j}^{n}} \varphi d\mu \\
&= \int_{\Rd} \varphi(x) K_{Z,M}(x,A) d\mu(x).
\end{aligned}.
\end{equation}
For a general $\phi \in \mathcal{M}_{B,c}(\R^{p})$, we use a sequence $(\phi_{n})_{n} \subset \mathcal{M}_{B,c}(\R^{p})$ of simple functions converging uniformly to $\phi$, say $\phi_{n} = \sum_{j \in \tilde{J}_{n}} a_{j}^{n} \mathbbm{1}_{A_{j}^{n}}$. Then,
\begin{equation}
\label{Eq:ProofIteratedIntKCrossPos}
\begin{aligned}
\Cov\left( \int_{\Rd}Z(x)\varphi(x)d\mu(x) , \int_{\R^{p}} \phi(y)dM(y) \right) &= \lim_{n \to \infty} \Cov\left( \int_{\Rd}Z(x)\varphi(x)d\mu(x) , \sum_{j \in \tilde{J}_{n}} a_{j}^{n} M(A_{j}^{n})   \right)  \\
&=\lim_{n \to \infty} \sum_{j \in \tilde{J}_{n}} \overline{a_{j}^{n}} \int_{\Rd}K_{Z,M}(x , A_{j}^{n})\varphi(x)d\mu(x) \\
&= \lim_{n \to \infty} \int_{\R^{p}} \overline{\phi_{n}}(y) d\left( \int_{\Rd} \varphi(x)K_{Z,M}(x,\cdot)d\mu(x) \right)(y)  \\
&= \int_{\R^{p}} \overline{\phi}(y) d\left( \int_{\Rd} \varphi(x)K_{Z,M}(x,\cdot)d\mu(x) \right)(y).
\end{aligned}
\end{equation}
The final result comes from Fubini Theorem \eqref{Theo:FubiniKernelMu}. Note that we have used measure \eqref{Eq:IntVarphiKmux(informal)}. $\blacksquare$

\section{The self-integral}
\label{Sec:SelfIntegral}

This Section is mainly devoted to introduce new concepts. We consider a cf-m kernel $K$ acting over $\Rd\times\Rd$ (the same space in both components). In such case we shall simply say that $K$ is a \textit{cf-m kernel over} $\Rd$. Let $D \in \BoundedBorelRd$. We would like to \textit{integrate $K$ with respect to itself}  over $D$, integrating its first component with respect to the second one. This number, if well-defined, will be written as
\begin{equation}
\label{Eq:SelfIntKinformal}
``\int_{D}K(x,dx)".
\end{equation}
A trivial example where this peculiar idea works is the following: consider $K$ of the form $K(x,A) = f(x)\mu(A)$, with $f \in C(\Rd)$ and $\mu \in \mathscr{M}(\Rd)$. In such case, we simply do
\begin{equation}
\label{Eq:SelfIntTensor}
\int_{D}K(x,dx) := \int_{D}fd\mu.
\end{equation}
The same idea can be applied if $K$ is a finite linear combination of kernels of this form $K = \sum_{j} f_{j} \otimes \mu_{j}$, that is, members of $C(D)\otimes\mathscr{M}(D)$, for which the self-integral can be defined unambiguously as
\begin{equation}
\label{Eq:SelfIntegralTensorLC}
\int_{D} K(x,dx) := \sum_{j} \int_{D} f_{j} d\mu_{j}.
\end{equation}
An analyst may have the reflex extending the definition to a completion of the space $C(D)\otimes\mathscr{M}(D)$ with a suitable topology. This approach will be considered in a forthcoming paper. For now, we shall use Riemann sums to try to define the self-integral. Let $\RiemannSystem{I}{x}{J}$   be a Riemann  system of $D$. Grounded on the definition of Riemann-Stieltjes integrals, we ask ourselves if the limit
\begin{equation}
\label{Eq:SelfIntegralDefRiemmanSums}
\int_{D}K(x,dx) := \lim_{n \to \infty} \sum_{j \in J_{n}}K(x_{j}^{n} , I_{j}^{n} ),
\end{equation}
exists and if it is independent of the Riemann system. One fact that many people know, but which may have not been expressed in this way, is that such limit is not always well-defined.

\begin{Prop}
\label{Prop:ExampleNotDefSelfInt}
There are cross-pos-def kernels for which the self-integral \eqref{Eq:SelfIntegralDefRiemmanSums} is not uniquely-defined.
\end{Prop} 

\textbf{Proof:} We give as counterexample the cross-covariance kernel between a White Noise and its primitive when in $d=1$ (Brownian motion)\footnote{This is the classical difference between Itô and Stratonovitch integrals}. Consider $D = [0,1)$ and the kernel $K(x,A) = \lambda( A \cap [0,x] )$. Then, if we use partitions of $[0,1)$ with $n$ same-length intervals $I_{j}^{n} = [ \frac{j-1}{n} , \frac{j}{n})$, $j \in J_{n} = \lbrace 1 , ... , n \rbrace $, and if we use the left-points of the intervals as tag-points, we obtain
\begin{equation}
\label{Eq:SelfIntegralIto}
\int_{[0,1)}K(x,dx) = \lim_{n \to \infty} \sum_{j \in J_{n}} \lambda\left( [\frac{j-1}{n} , \frac{j}{n})\cap [0 , \frac{j-1}{n} ] \right) = 0.
\end{equation} 
But if we use the middle-points as tag-points, we obtain
\small
\begin{equation}
\label{Eq:SelfIntegralStrat}
\int_{[0,1)} K(x,dx) =  \lim_{n \to \infty} \sum_{j \in J_{n}} \lambda\left( [\frac{j-1}{n} , \frac{j}{n})\cap [0 , \frac{2j-1}{2n} ] \right) =   \lim_{n \to \infty} \sum_{j \in J_{n}} \lambda\left( [\frac{j-1}{n} , \frac{2j-1}{2n} ] \right) = \frac{1}{2}.
\end{equation}
\normalsize
Hence, the limit depends upon the Riemann system employed. $\blacksquare$ 

The previous result inspires the next definition.

\begin{Def}
\label{Def:SelfIntegral}
$K$ is said to be self-integrable over $D$ if there exists a unique $\ell_{D} \in \R$ such that for every Riemann system of $D$, say $\RiemannPartitionTag{I}{x}{J}$, it holds
\begin{equation}
\label{Eq:DefSelfIntegralK}
\ell_{D} = \lim_{n \to \infty} \sum_{j \in J_{n}} K(x_{j}^{n} , I_{j}^{n} ).
\end{equation}
In such case, $\ell_{D}$ is called the self-integral of $K$ over $D$, and it is denoted
\begin{equation}
\label{Eq:SelfIntegralNotation}
\int_{D}K(x,dx) := \ell_{D}.
\end{equation}
\end{Def}

Plenty of questions may be asked concerning this concept, but we shall only mention some of them in the conclusive Section \ref{Sec:OpenQuestions}, since in this work we shall focus on the application of this concept in stochastic integration. In order to contrast with the not-self-integrable situation of Proposition \ref{Prop:ExampleNotDefSelfInt}, we give the following non-trivial example of well-defined self-integral.

\begin{Exa}
\label{Ex:SelfIntegralIteratedKintegrals}
Consider the weirdly-noted integral \eqref{Eq:IntPsiKxy(pd)(informal)} for which we gave a rigorous meaning through the use of the measure $\tilde{K}_{\mu}$ obtained in Proposition \ref{Prop:ExistenceMeasureKmu}. We will show that it can be conceived as a self-integral. Consider the case $\psi \in C_{c}(\Rd\times\R^{p})$ with $\supp(\psi) \subset D_{d}\times D_{p}$. Then, we can consider the kernel $K_{\psi,\mu} : \R^{p} \times \mathcal{B}_{B}(\R^{p}) \to \R$ given by 
\begin{equation}
\label{Eq:KernelPsiMuSelfIntegral}
K_{\psi,\mu}(y,A) = \int_{\Rd}\psi(x,y)K(x,A)d\mu(x),
\end{equation}
which is a cf-m kernel if $K$ is locally bounded. Note that $K_{\psi,\mu} = 0 $ for $x \notin D_{p}$. We claim that $K_{\psi,\mu}$ is self-integrable over $D_{p}$, having
\begin{equation}
\label{Eq:SelfIntegralKPsiMu}
\int_{D_{p}} K_{\psi,\mu}(y,dy) = \int_{\R^{p}} d\left( \int_{\Rd} \psi(x,\cdot)K(x,\cdot) d\mu(x) \right).
\end{equation}
Indeed, if $\RiemannPartitionTag{I}{y}{J}$ is any Riemann system of $D_{p}$, then
\begin{equation}
\begin{aligned}
\lim_{n \to \infty} \sum_{j \in J_{n}} K_{\psi , \mu}(y_{j}^{n} , I_{j}^{n} ) &= \lim_{n \to \infty} \sum_{j \in J_{n}} \int_{\Rd} \psi( x , y_{j}^{n} )K( x , I_{j}^{n} ) d\mu(x) \\
&= \int_{\Rd} \lim_{n \to \infty} \sum_{j \in J_{n}} \psi( x , y_{j}^{n} )K( x , I_{j}^{n} ) d\mu(x) \\
&= \int_{\Rd} \int_{\R^{p}} \psi(x,y)dK(x,\cdot)(y) d\mu(x) = \int_{\R^{p}} d\left( \int_{\Rd} \psi(x,\cdot)K(x,\cdot) d\mu(x) \right),
\end{aligned}
\end{equation}
where we have used LCDT and Fubini Theorem \ref{Theo:FubiniKernelMu}. The Riemann system of $D_{p}$ was arbitrary, thus the kernel $K_{\psi,\mu}$ is self-integrable over $D_{p}$ and the self-integral coincides with the iterated integral \eqref{Eq:IntPsiKxy(pd)(informal)}. $\square$ 
\end{Exa}

We shall now introduce an auxiliary concept which plays an important role in future developments.

\begin{Def}
\label{Def:SecondOrderKernel}
Let $K$ be a cf-m kernel over $\Rd\times\R^{p}$. The second-order kernel $K^{(2)}$ is defined as the kernel over $(\Rd\times\Rd)\times(\R^{p}\times\R^{p})$ given by $K^{(2)}((x,y),E ) := ( K(y,\cdot)\otimes K(x,\cdot) )(E)$.
\end{Def}

The second-order kernel is a kind of cross-tensor product measure (remark the order $(y,x)$ and not $(x,y)$ in Definition \ref{Def:SecondOrderKernel}). It is clearly well-defined as a function-measure kernel, though it is not clear if it is in general continuous in $(x,y)$. We shall see later that this is true for cross-pos-def kernels, which is the case we shall focus on in the following sections. For now, we shall introduce the next Definition.

\begin{Def}
\label{Def:QuasiSelfIntegral}
Suppose $K$ is such that $K^{(2)}$ is a cf-m kernel. Let $A,B \in \BoundedBorel{\Rd}$. We say that $K^{(2)}$ is quasi-self-integrable over $A\times B$ if there exists $\ell_{A,B} \in \R$ such that for every Riemann system of $A$, $\RiemannPartitionTag{A}{a}{J^{A}}$, and every Riemann system of $B$, $\RiemannPartitionTag{B}{b}{J^{B}}$, the following double-limit holds
\begin{equation}
\label{Eq:DefDoubleLimitQuasiSelfIntegral}
\lim_{(n,m)\to \infty} \sum_{j \in J_{n}}\sum_{k \in J^{B}_{m}} K^{(2)}( (a_{j}^{n} , b_{k}^{m} ) , A_{j}^{n} \times B_{k}^{m}  ) = \ell_{A,B}.
\end{equation}
In such case, we call $\ell_{A,B}$ the quasi-self-integral of $K^{(2)}$ over $A\times B$, and we denote it
\begin{equation}
\label{Eq:QuasiSelfIntegralNotation}
\tilde{\int}_{A\times B} K^{(2)}(x,dx) := \ell_{A,B}.
\end{equation}
\end{Def}

Definition \ref{Def:QuasiSelfIntegral} is weaker than the definition of a proper self-integral of $K^{(2)}$ since it only considers a specific class of Riemann system: partitions constructed with products $A_{j}^{n}\times B_{k}^{m}$. In addition, Definition \ref{Def:QuasiSelfIntegral} only considers integration over sets of the form $A\times B$, rather than general bounded Borel sets of $\RdxRd$. It follows that if $K^{(2)}$ is self-integrable over $E \in \BoundedBorel{\Rd\times\R^{p}}$, with $E = A\times B$, then it is quasi-self-integrable over $E$. In such case, the quasi-self-integral and the self-integral coincide. It is not clear if quasi-self-integrability implies self-integrability.

\section{Self-integrability and stochastic integration}
\label{Sec:SelfIntegrabilityAndStochasticIntegration}

In this Section we show that the concept of self-integrability plays a central role in the definition of stochastic integrals. All along this Section, $Z$ is a centred mean-square continuous stochastic process over $\Rd$, $M$ is a centred random measure over $\Rd$, and $D \in \BoundedBorel{\Rd}$.

\subsection{A first approach}
\label{Sec:FirstApproachStochIntegrals}

The first simple question is whether we can define the stochastic integral
\begin{equation}
``\int_{D}Z(x)dM(x)".
\end{equation}
Analogously to the deterministic case, we study the limit of the  random variables
\begin{equation}
\label{Eq:LimRiemannSumsZM}
\lim_{n \to \infty} \sum_{j \in J_{n}} Z(x_{j}^{n})M(I_{j}^{n}), 
\end{equation}
being $\RiemannPartitionTag{I}{x}{J}$ a Riemann system of $D$. We introduce hence the following definition.

\begin{Def}
\label{Def:StochasticIntegralUnique}
We say that the stochastic integral of $Z$ with respect to $M$ over $D$ is uniquely-defined if the limit \eqref{Eq:LimRiemannSumsZM} exists in  the sense of $L^{p}\OmegaAP$ for a specified $p \in [1,\infty]$ and that it does not depend upon the Riemann system of $D$ used. In such case we denote
\begin{equation}
\label{Eq:DefStochasticIntegral}
\int_{D}ZdM := \int_{D}Z(x)dM(x) := \lim_{n \to \infty} \sum_{j \in J_{n}} Z(x_{j}^{n})M(I_{j}^{n}).
\end{equation}
\end{Def}

In this work we will use $p=1 $ or $p=2$. Note that the random variables in \eqref{Eq:LimRiemannSumsZM} are always in $\LpOmegaAP{1}$, hence it does make sense to ask if they converge in the $\LpOmegaAP{1}$ sense. However, in some cases the convergence in $\LpOmegaAP{2}$ is easier to analyse (such convergence implies of course the $\LpOmegaAP{1}$ convergence). We shall always specify the sense of the convergence we are considering. 

The first link between self-integrability and stochastic integrals is given by the next Proposition.

\begin{Prop}
\label{Prop:SelfIntegralNecessary}
If $\int_{D}ZdM$ is uniquely-defined in the $\LpOmegaAP{p}$ sense for any $p \in [1,\infty]$, then the cross-covariance kernel $K_{Z,M}$ is self-integrable over $D$.
\end{Prop}

\textbf{Proof:} If the limit \eqref{Eq:LimRiemannSumsZM} converges to some random variable in $L^{p}\OmegaAP$, then the expectations also converge to some real limit $\ell \in \R$. Hence
\begin{equation}
\label{Eq:MeanSelfIntegralConvergence}
\lim_{n \to \infty}  \mathbb{E}\left( \sum_{j \in J_{n}} Z(x_{j}^{n} ) M( I_{j}^{n} ) \right) = \lim_{n \to \infty} \sum_{j \in J_{n}} K_{Z,M}(x_{j}^{n},I_{j}^{n} ) = \ell.
\end{equation}
Since this holds for every Riemann system of $D$, $K_{Z,M}$ must be self-integrable over $D$. $\blacksquare$

The last Proposition does not only tell us that we need self-integrability of $K_{Z,M}$ in order to define uniquely the stochastic integral, but it also gives us the expectation formula
\begin{equation}
\label{Eq:MeanSelfIntegral}
\mathbb{E}\left( \int_{D}ZdM \right) = \int_{D}K_{Z,M}(x,dx).
\end{equation}
There are simple cases where the stochastic integral is uniquely-defined and hence the associated cross-covariance kernel is self-integrable. We list three important examples which will be used further.

\begin{Exa}
\label{Ex:StochasticIntegralZMindep}
Consider $Z$ and $M$ independent. Let us prove the Riemann sums \eqref{Eq:LimRiemannSumsZM} form  a Cauchy sequence in $\LtwoOmega$. Using the measure structure of $M$, we right the Cauchy gaps as
\begin{equation}
\label{Eq:CauchyGapStochRiemannSum}
\sum_{j \in J_{n}} Z(x_{j}^{n})M(I_{j}^{n} ) - \sum_{k \in J_{m}} Z(x_{k}^{m})M(I_{k}^{m} ) = \sum_{j \in J_{n}} \sum_{k \in J_{m}} ( Z(x_{j}^{n}) - Z(x_{k}^{m}) )M(I_{j}^{n} \cap I_{k}^{m} ).
\end{equation}
Taking $\mathbb{E}\left[ ( \cdot )^{2} \right]$ of this, one obtains the multi-sum
\begin{equation}
\label{Eq:2MomentCauchyGapStochRiemannSum}
\sum_{j,j' \in J_{n}} \sum_{k,k' \in J_{m}} \mathbb{E}\left[  ( Z(x_{j}^{n}) - Z(x_{k}^{m}) )( Z(x_{j'}^{n}) - Z(x_{k'}^{m}) )M(I_{j}^{n} \cap I_{k}^{m} )M(I_{j'}^{n} \cap I_{k'}^{m} ) \right].
\end{equation}
From the independence condition and using the increment process $I_{Z}$ (See \eqref{Eq:IncrementProcess} and \eqref{Eq:CovIncrementProcess}), we obtain
\begin{equation}
\label{Eq:MultiSumCauchyGapZMindep}
\sum_{j,j' \in J_{n}} \sum_{k,k' \in J_{m}} C_{I_{Z}}( (x_{j}^{n} , x_{k}^{m}),(x_{j'}^{n} , x_{k'}^{m})  ) C_{M}( ( I_{j}^{n} \cap I_{k}^{m}) \times  (I_{j'}^{n} \cap I_{k'}^{m})  ).
\end{equation}
Now we consider the uniform continuity of $C_{I_{Z}}$ over $\overline{D}^{4}$. Note that $C_{I_{Z}}( (x,x) , (y,y) ) = 0$. Hence, for every $\epsilon > 0 $ there exists $\delta > 0 $ such that if $|x-u| + |y-v|< \delta$, then $|C_{I_{Z}}( (x,u),(y,v) )| < \epsilon $. Consider hence $n_{0}$ such that if $n,m \geq n_{0}$, $\diam( I_{j}^{n} ) + \diam(I_{k}^{m}) < \frac{\delta}{2} $ for every $(j,k) \in J_{n}\times J_{m}$. Thus, we claim that if $n,m \geq n_{0}$ then for every $j,j \in J_{n}$ and every $k,k'\in J_{m}$ one has
\begin{equation}
\label{Eq:BoundindCIzCM}
|C_{I_{Z}}( (x_{j}^{n} , x_{k}^{m}),(x_{j'}^{n} , x_{k'}^{m})  ) C_{M}( ( I_{j}^{n} \cap I_{k}^{m}) \times  (I_{j'}^{n} \cap I_{k'}^{m})  )| < \epsilon |C_{M}|( (I_{j}^{n}\times I_{k}^{m}) \cap (I_{j}^{n}\times I_{k}^{m}) ).
\end{equation}
This holds since if $I_{j}^{n}\cap I_{k}^{m} = \emptyset $ then \eqref{Eq:BoundindCIzCM} equals $0$, and if $I_{j}^{n}\cap I_{k}^{m} \neq \emptyset$, using  $\diam( I_{j}^{n}) + \diam(I_{k}^{m}) < \frac{\delta}{2}$, we conclude $|x_{j}^{n} - x_{k}^{m}| < \frac{\delta}{2}$, and the same idea can be used for $|x_{j'}^{n}-x_{k'}^{m}|$. Hence $|C_{I_{Z}}( (x_{j}^{n} , x_{k}^{m}),(x_{j'}^{n} , x_{k'}^{m})  )| < \epsilon$. It follows that \eqref{Eq:MultiSumCauchyGapZMindep} is bounded by
\begin{equation}
\sum_{j,j' \in J_{n}} \sum_{k,k' \in J_{m}} \epsilon |C_{M}|( ( I_{j}^{n} \cap I_{k}^{m}) \times  (I_{j'}^{n} \cap I_{k'}^{m})  ) = \epsilon |C_{M}|(D\times D).
\end{equation}
The Riemann sums form thus a Cauchy sequence in $\LtwoOmega$, converging to a random variable with finite variance. If $\RiemannSystem{D}{y}{\tilde{J}}$ is another Riemann system of $D$, we can consider the gap
\begin{equation}
\label{Eq:GapRiemannSumsDiffSystem}
\sum_{j \in J_{n}} Z(x_{j}^{n})M(I_{j}^{n} ) - \sum_{k \in \tilde{J}_{n}} Z(y_{k}^{n})M(D_{k}^{n} ) = \sum_{j \in J_{n}} \sum_{k \in \tilde{J}_{n}} ( Z(x_{j}^{n}) - Z(y_{k}^{n}) )M(I_{j}^{n} \cap D_{k}^{n} ).
\end{equation}
An analogue argument can be used to prove that the $\LtwoOmega$-norm of \eqref{Eq:GapRiemannSumsDiffSystem} goes to $0$ as $n \to \infty$. The limit is hence independent of the Riemann system chosen, so the stochastic integral $\int_{D}ZdM$ is uniquely-defined in the sense of $\LtwoOmega$. Of course in this case $K_{Z,M} = 0$, which is trivially self-integrable. $\square$
\end{Exa}

\begin{Exa}
\label{Ex:FiniteImage}
Suppose $Z$ or $M$ has finite-dimensional image in $\LtwoOmega$. For instance, suppose $Z$ is of the form
\begin{equation}
\label{Eq:ZfiniteSum}
Z = \sum_{\alpha \leq m} X_{\alpha}f_{\alpha},
\end{equation}
where $(X_{\alpha})_{\alpha \leq m} \subset \LtwoOmega$ and $(f_{\alpha})_{\alpha \leq m} \subset C(\Rd)$. Then, the stochastic integral $\int_{D}ZdM$ is uniquely-defined since for every Riemann system $\RiemannSystem{I}{x}{J}$ we have
\begin{equation}
\label{Eq:IntZdMFiniteSum}
\int_{D}ZdM = \lim_{n \to \infty} \sum_{j \in J_{n}} Z(x_{j}^{n})M(I_{j}^{n}) = \sum_{\alpha \leq m} X_{\alpha} \lim_{n \to \infty} \sum_{j \in J_{n}} f_{\alpha}(x_{j}^{n}) M(I_{j}^{n} ) =  \sum_{\alpha \leq m} X_{\alpha} \int_{D}f_{\alpha}(x)dM(x).
\end{equation}
In this case the cross-covariance kernel $K_{Z,M}$ has a simple form. Let $\sigma_{X_\alpha}^{2} = \Var(X_{\alpha}) > 0$, and let $\nu_{\alpha}(A) :=  \sigma_{X_{\alpha}}^{-1}\mathbb{E}( X_{\alpha} M(A) )$. Form the $\sigma$-additivity of $M$, $\nu_{\alpha}$ defines a measure over $\Rd$. Thus,
\begin{equation}
\label{Eq:CrossCovarianceFiniteSum}
K_{Z,M}(x,A) = \mathbb{E}( Z(x)M(A)) = \sum_{\alpha \leq m} \sigma_{X_{\alpha}} f_{\alpha}(x) \nu_{\alpha}(A).
\end{equation}
$K_{Z,M}$ is hence a linear combination of tensors. It is then self-integrable over $D$. Analogously, this same situation is present when the random measure $M$ is such that
\begin{equation}
\label{Eq:MFiniteSum}
M = \sum_{\beta \leq m} Y_{\beta} \mu_{\beta},
\end{equation}
with $(Y_{\beta})_{\beta \leq m} \in \LtwoOmega$, and $(\mu_{\beta})_{\beta \leq m} \in \mathscr{M}(\Rd)$. We shall explode more in detail this kind of construction in a forthcoming paper. $\square$
\end{Exa}

\begin{Exa}
\label{Ex:MRegular}
Consider the case where $M$ is absolutely continuous with respect to a real deterministic measure $\mu \in \mathscr{M}(\Rd)$ in the sense that it exists an enough-regular process ($\mu$-almost everywhere mean-square continuous-locally-bounded, for instance) $U$ such that
\begin{equation}
\label{Eq:MregularIntUmu}
M(A) = \int_{A}U(y)d\mu(y), \quad \forall A \in \BoundedBorelRd.
\end{equation}
The Riemann sums approaching the stochastic integrals are given by
\begin{equation}
\label{Eq:RiemannSum}
\sum_{j \in J_{n}} Z(x_{j}^{n})\int_{I_{j}^{n}} U(y) d\mu(y) 
\end{equation}
When analysing the $\LpOmegaAP{1}$-norm of the Cauchy gap \eqref{Eq:CauchyGapStochRiemannSum} we obtain
\small
\begin{equation}
\label{Eq:L1CauchyGapintZUdmu}
\begin{aligned}
\mathbb{E}\Big( \Big|  \sum_{j \in J_{n}} \sum_{k \in J_{m}} (Z(x_{j}^{n}&) - Z(x_{k}^{m}) )\int_{I_{j}^{n}\cap I_{k}^{m}}U(y)d\mu(y)  \Big|  \Big) \leq \sum_{j \in J_{n}} \sum_{k \in J_{m}} \mathbb{E}\Big( |Z(x_{j}^{n}) - Z(x_{k}^{m}) | \Big|  \int_{I_{j}^{n}\cap I_{k}^{m}}U(y)d\mu(y)   \Big|  \Big) \\
&\leq \sum_{j \in J_{n}} \sum_{k \in J_{m}} \sqrt{\mathbb{E}( |Z(x_{j}^{n}) - Z(x_{k}^{m}) |^{2} )} \sqrt{\mathbb{E}\Big( \Big|  \int_{I_{j}^{n}\cap I_{k}^{m}}U(y)d\mu(y)   \Big|^{2}  \Big)   } \\
&= \sum_{j \in J_{n}} \sum_{k \in J_{m}} \sqrt{C_{I_{Z}}(  (x_{j}^{n} , x_{k}^{m} ) , (x_{j}^{n} , x_{k}^{m} )  )}  \sqrt{\int_{I_{j}^{n}\cap I_{k}^{m}} \int_{I_{j}^{n}\cap I_{k}^{m} } C_{U}(y,z) d\mu(y) d\mu(z)  },
\end{aligned}
\end{equation}
\normalsize
where $C_{U}$ is the covariance function of $U$. We can use the Cauchy-Schwarz inequality for $C_{U}$ to obtain
\begin{equation}
\begin{aligned}
\int_{I_{j}^{n}\cap I_{k}^{m}} \int_{I_{j}^{n}\cap I_{k}^{m} } C_{U}(y,z) d\mu(y) d\mu(z)  &\leq  \int_{I_{j}^{n}\cap I_{k}^{m}} \int_{I_{j}^{n}\cap I_{k}^{m} } \sqrt{C_{U}(y,y)}\sqrt{C_{U}(z,z)} d|\mu|(y) d|\mu|(z) \\
&= \left( \int_{I_{j}^{n} \cap I_{k}^{m}} \sqrt{C_{U}(y,y)}d|\mu|(y) \right)^{2}.
\end{aligned}
\end{equation}
With the same arguments as in \eqref{Eq:BoundindCIzCM}, for any $\epsilon > 0 $ and for $n,m$ big enough \eqref{Eq:L1CauchyGapintZUdmu} is bounded by
\begin{equation}
\sum_{j \in J_{n}} \sum_{k \in J_{m}} \epsilon\int_{I_{j}^{n} \cap I_{k}^{m}} \sqrt{C_{U}(y,y)}d|\mu|(y) = \epsilon \int_{D} \sqrt{C_{U}(y,y)} d|\mu|(y).  
\end{equation}
Thus the Riemann sums form a Cauchy sequence on $\LpOmegaAP{1}$ hence they converge. A similar algebraic computation assures that the limit does not depend upon the Riemann system chosen. The limit in the sense of $\LpOmegaAP{1}$ is hence $\int_{D}ZdM$. In this case the cross-covariance kernel $K_{Z,M}$ has the form
\begin{equation}
K_{Z,M}(x,A) = \int_{A} C_{Z,U}(x,y)d\mu(y),
\end{equation}
where $C_{Z,U}$ is the \textit{cross-covariance function} between $Z$ and $U$. The self-integral of $K_{Z,M}$ is given by
\begin{equation}
\int_{D}K_{Z,M}(x,dx) = \int_{D}C_{Z,U}(x,x)d\mu(x). 
\end{equation}
Note that with this analysis we can define $\int_{D} Z(x)U(x) d\mu(x) := \int_{D}ZdM$, which should coincide with other intuitive manners of constructing such integral.$\square$
\end{Exa}

\subsection{Gaussian case: necessary and sufficient conditions}
\label{Sec:GaussianCase}

Consider the case where $Z$ and $M$ are jointly Gaussian. That is, for every $(x_{1},... ,x_{n}) \in (\Rd)^{n}$ and every $(A_{1}, ... , A_{m}) \in (\BoundedBorel{\Rd})^{m}$, the vector $(Z(x_{1}) , ... , Z(x_{n}) , M(A_{1}) , ... , M(A_{m}) )$ is a Gaussian vector. In such case,  $C_{Z}, C_{M}$ and $K_{Z,M}$ describe entirely the finite-dimensional laws of the random vectors we can construct with $Z$ and $M$. Concerning the unique definition of $\int_{D}ZdM$, note that in Examples \ref{Ex:FiniteImage} and \ref{Ex:MRegular} there are extra regularity conditions on $Z$ or $M$ playing a role. In Example \ref{Ex:StochasticIntegralZMindep} there is no regularity requirement on $Z$ or $M$, but a strong independence condition is used. It is then expected that conditions for the definition of $\int_{D}ZdM$, without asking more regularity in $Z$ or $M$, will rely entirely on the \textit{dependence} between $Z$ and $M$, which is fully determined by $K_{Z,M}$. The main result of this work is Theorem \ref{Theo:IntZMGaussian} which shows that this intuition is correct. For introducing such Theorem we will need Definition \ref{Def:QuasiSelfIntegral} applied to the kernel $K_{Z,M}$. The next Proposition will helps us for that.
\begin{Prop}
\label{Prop:KsecondorderCrossPosDef}
If $K$ is cross-pos-def, then $K^{(2)}$ is also cross-pos-def.
\end{Prop}
Proposition \ref{Prop:KsecondorderCrossPosDef} and Theorem \ref{Theo:KcrossPosDefEquicontinuous} imply that $K_{Z,M}^{(2)}$ is a cf-m kernel, hence Definition \ref{Def:QuasiSelfIntegral} of the quasi-self-integral can be applied to it. In order to prove Proposition \ref{Prop:KsecondorderCrossPosDef} we shall use the following Lemma which introduces the tensor product of independent random measures.

\begin{Lemma}
\label{Lemma:M1tensorM2indep}
Let $M_{1}$ and $M_{2}$ be two independent random measures over $\Rd$ and $\R^{p}$ respectively. Then, there exists a unique random measure $M_{1}\otimes M_{2}$ over $\Rd\times\R^{p}$ satisfying
\begin{equation}
\label{Eq:DefM1OtimesM2indep}
(M_{1}\otimes M_{2})(A\times B) = M_{1}(A)M_{2}(B), \quad \forall A \in \BoundedBorelRd, B \in \mathcal{B}_{B}(\R^{p}).
\end{equation}
\end{Lemma}

\textbf{Proof of Lemma \ref{Lemma:M1tensorM2indep}:} We define the random measure $M_{1}\otimes M_{2}$ acting over functions in $C_{c}(\Rd\times\R^{p})$ through the iterated stochastic integrals
\begin{equation}
\label{Eq:DefM1tensorM2IteratedInt}
\langle M_{1}\otimes M_{2} , \psi \rangle := \int_{\Rd} \int_{\R^{p}} \psi(x,y) dM_{2}(y) dM_{1}(x), \quad \forall \psi \in C_{c}(\Rd\times\R^{p}).
\end{equation}
Note that the stochastic process $ x\mapsto \int_{\R^{p}} \psi(x,y) dM_{2}(y)$ is mean-square continuous and independent of $M_{1}$, and hence we are in the case of Example \ref{Ex:StochasticIntegralZMindep} of a uniquely-defined stochastic integral. The application $\psi \mapsto \langle M_{1}\otimes M_{2} , \psi \rangle$ has null expectation (Eq. \eqref{Eq:MeanSelfIntegral}) and covariance
\begin{equation}
\label{Eq:CovM1TensorM2Indep}
\Cov(\langle M_{1}\otimes M_{2} , \psi_{1} \rangle , \langle M_{1}\otimes M_{2}  , \psi_{2} \rangle  ) = \int_{\Rd\times \Rd} \int_{\R^{p}\times\R^{p}} \psi_{1}(x,u)\overline{\psi_{2}}(y,v)  dC_{M_{2}}(u,v) dC_{M_{1}}(x,y), 
\end{equation}
which defines a measure over $(\Rd\times\R^{p})\times(\Rd\times\R^{p})$ since $C_{M_{1}}$ and $C_{M_{2}}$ are measures. From Theorem \ref{Theo:ExtensionRandomMeasure}, $M_{1}\otimes M_{2}$ defines a unique random measure over $\Rd\times\R^{p}$ satisfying \eqref{Eq:DefM1OtimesM2indep} (use $\psi = \mathbbm{1}_{A}\otimes \mathbbm{1}_{B}$). $\blacksquare$

\textbf{Proof of Proposition \ref{Prop:KsecondorderCrossPosDef}:} Let $Z_{1},Z_{2}$ be centred mean-square continuous processes over $\Rd$ and let $M_{1}$ and $M_{2}$ be centred random measures over $\R^{p}$ such that $(Z_{1},M_{1})$ is independent of $(Z_{2},M_{2})$ and such that $K$ is the cross-covariance kernel of both pairs, $K = K_{Z_{1},M_{1}} = K_{Z_{2},M_{2}}$. Consider thus the centred and mean-square continuous process $U(x,y) := Z_{2}(x)Z_{1}(y)$ and the random measure $M_{1}\otimes M_{2}$. Then,
\begin{equation}
\begin{aligned}
\mathbb{E}\big(U(x,y)(M_{1}\otimes M_{2})(A\times B) \big) &= \mathbb{E}\big( Z_{1}(y)M_{1}(A)  \big) \mathbb{E}\big( Z_{2}(x)M_{2}(B)  \big) \\
&= K(y,A)K(x,B) \\
&= K^{(2)}((x,y) , A\times B ),
\end{aligned}
\end{equation}
for every $A \in \BoundedBorelRd$, $B \in \mathcal{B}_{B}(\R^{p})$. Thus, $K^{(2)}$ is nothing but the cross-covariance kernel between  $U$ and $M_{1}\otimes M_{2}$, hence it is cross-pos-def. $\blacksquare$

We state now Theorem \ref{Theo:IntZMGaussian}.

\begin{Theo}
\label{Theo:IntZMGaussian}
If $Z$ and $M$ are jointly Gaussian, then $\int_{D}ZdM$ is uniquely-defined in a $\LtwoOmega$ sense if and only if $K_{Z,M}$ is self-integrable over $D$ and $K^{(2)}_{Z,M}$ is quasi-self-integrable over $D\times D$.
\end{Theo}

We shall use next two simple Lemmas, proper of basic Probability Theory.

\begin{Lemma}
\label{Lemma:Cov4MomentGaussian}
If $(X_{1},X_{2},X_{3},X_{4})$ is a centred Gaussian vector, then 
\begin{equation}
\Cov( X_{1}X_{2} , X_{3}X_{4}) = \mathbb{E}( X_{1}X_{3} )\mathbb{E}( X_{2}X_{4}) + \mathbb{E}( X_{1}X_{4} )\mathbb{E}( X_{2}X_{3}).
\end{equation}
\end{Lemma}

\textbf{Proof of Lemma \ref{Lemma:Cov4MomentGaussian}:} This is a particular case of Isserlis-Wick Theorem \citep{isserlis1918formula,wick1950evaluation}. $\blacksquare$

\begin{Lemma}
\label{Lemma:L2Convergence}
Let $(X_{n})_{n}, (Y_{n})_{n}$ be two sequences in $\LpOmegaAP{2}$. Then,
\begin{enumerate}[(i)]
\item \label{It:L2Convergence1}$(X_{n})_{n}$ converges in $\LpOmegaAP{2}$ if and only if $\displaystyle\lim_{n \to \infty} \mathbb{E}(X_{n})$ and $\displaystyle\lim_{(n,m) \to \infty} \Cov(X_{n},X_{m})$ exist.
\item \label{It:L2COnvergence2}If $X_{n} \stackrel{L^{2}(\Omega)}{\to} X $, then $Y_{n} \stackrel{L^{2}(\Omega)}{\to} X $ if and only if  $\displaystyle\lim_{n \to \infty} \mathbb{E}(Y_{n}) = \displaystyle\lim_{n \to \infty} \mathbb{E}(X_{n})$ and $\displaystyle\lim_{(n,m) \to \infty} \Cov(Y_{n},Y_{m})  = \displaystyle\lim_{(n,m) \to \infty} \Cov(X_{n},X_{m}) = \displaystyle\lim_{(n,m) \to \infty} \Cov(Y_{n},X_{m})$.
\end{enumerate} 
\end{Lemma}

\textbf{Proof of Lemma \ref{Lemma:L2Convergence}:} For $\ref{It:L2Convergence1}$, the necessity is clear in the case of the mean, and for the covariance it is enough to consider $\Cov(X_{n},X_{m}) = \frac{1}{2}\left(\Var(X_{n}) + \Var(X_{m}) - \Var(X_{n}-X_{m}) \right)$, which converges to the variance of the limit as $(n,m) \to \infty$. For the sufficiency, the convergence of the double-limit of the covariances implies in particular (using $n=m$) the convergence of the variances $\Var(X_{n})$, whose limit coincide with $\lim_{(n,m) \to \infty} \Cov( X_{n},X_{m})$. Using $\mathbb{E}(|X_{n}-X_{m}|^{2}) = \Var(X_{n}) + \Var(X_{m}) - 2 \Cov(X_{n} , X_{m}) + |\mathbb{E}(X_{n}) - \mathbb{E}(X_{m}) |^{2}$, it follows that $(X_{n})_{n}$ is a Cauchy sequence hence convergent.

For $\ref{It:L2COnvergence2}$ the necessity in the case of the means is clear as well as the equality $\displaystyle\lim_{(n,m) \to \infty} \Cov(Y_{n},Y_{m})  = \displaystyle\lim_{(n,m) \to \infty} \Cov(X_{n},X_{m})$ (the limit is $\Var(X)$), and the equality involving $\Cov(Y_{n},X_{m})$ follows using  $\displaystyle\lim_{(n,m) \to \infty}\mathbb{E}(|X_{m}-Y_{n}|^{2}) = 0$ and developing. For the sufficiency, from $\ref{It:L2Convergence1}$ $(Y_{n})_{n}$ is convergent, and developing $\mathbb{E}(|Y_{n} - X_{n}|^{2})$ one proves that it goes to $0$ as $n \to \infty$, hence $(Y_{n})_{n}$ and $(X_{n})_{n}$ have the same limit. Details are left to the reader. $\blacksquare$

\textbf{Proof of Theorem \ref{Theo:IntZMGaussian}:} Let $\RiemannSystem{I}{x}{J}$, $\RiemannSystem{E}{y}{\tilde{J}}$ be two Riemann systems of $D$. The covariance between the Riemann sums approaching $\int_{D}ZdM$ is (Lemma \ref{Lemma:Cov4MomentGaussian})
\small
\begin{equation}
\label{Eq:CovDoubleRiemannSum}
\begin{aligned}
\Cov\Big( &\sum_{j \in J_{n}} Z(x_{j}^{n})M(I_{j}^{n}) , \sum_{k \in \tilde{J}_{m}} Z(y_{k}^{m})M(E_{k}^{m}) \Big) \\
&= \sum_{j \in J_{n}} \sum_{k \in \tilde{J}_{m}} \mathbb{E}( Z(x_{j}^{n})Z(y_{k}^{m}) ) \mathbb{E}(M(I_{j}^{n})M(E_{k}^{m})) + \mathbb{E}(Z(x_{j}^{n})M(E_{k}^{m}))\mathbb{E}(Z(y_{k}^{m})M(I_{j}^{n})) \\
&=  \sum_{j \in J_{n}} \sum_{k \in \tilde{J}_{m}} C_{Z}(x_{j}^{n} , y_{k}^{m}) C_{M}(I_{j}^{n}\times E_{k}^{m}) + K_{Z,M}(x_{j}^{n} , E_{k}^{m} ) K_{Z,M}(y_{k}^{m} , I_{j}^{n} ) \\
&= \sum_{j \in J_{n}} \sum_{k \in \tilde{J}_{m}} C_{Z}(x_{j}^{n} , y_{k}^{m}) C_{M}(I_{j}^{n}\times E_{k}^{m})  + \sum_{j \in J_{n}} \sum_{k \in \tilde{J}_{m}} K_{Z,M}^{(2)}( (x_{j}^{n} , y_{k}^{m} ), I_{j}^{n}\times E_{k}^{m} ).
\end{aligned}
\end{equation}
\normalsize
Note now that
\begin{equation}
\label{Eq:DoubleLimitIntCzCm}
\lim_{(n,m) \to \infty} \sum_{j \in J_{n}} \sum_{k \in \tilde{J}_{m}} C_{Z}(x_{j}^{n} , y_{k}^{m}) C_{M}(I_{j}^{n}\times E_{k}^{m}) = \int_{D\times D} C_{Z}dC_{M},
\end{equation}
which holds since the concerned double-sum is nothing but a Riemann approximation of the integral $\int_{D\times D} C_{Z}dC_{M}$.\footnote{If the reader needs more details, one can use the uniform continuity of $C_{Z}$ over $\overline{D}\times\overline{D}$ to prove that the function $\sum_{j \in J_{n}} \sum_{k \in \tilde{J}_{m}} C_{Z}(x_{j}^{n} , y_{k}^{m}) \mathbbm{1}_{I_{j}^{n}\times E_{k}^{m}} $ converges uniformly to $C_{Z}$ over $D\times D$ as $(n,m) \to \infty$ and then conclude.} It follows that
\small
\begin{equation}
\label{Eq:DoubleLimitCovRiemannSums}
\begin{aligned}
\lim_{(n,m) \to \infty} \Cov\Big( \sum_{j \in J_{n}} Z(x_{j}^{n})&M(I_{j}^{n}) \ , \sum_{k \in \tilde{J}_{m}} Z(y_{k}^{m})M(E_{k}^{m}) \Big) \\
&= \int_{D\times D} C_{Z}dC_{M} + \lim_{(n,m) \to \infty} \sum_{j \in J_{n}} \sum_{k \in \tilde{J}_{m}}  K_{Z,M}^{(2)}( (x_{j}^{n} , y_{k}^{m} ), I_{j}^{n}\times E_{k}^{m} ).
\end{aligned}
\end{equation}
\normalsize
We conclude hence the equivalence of the convergences of both double-limits in \eqref{Eq:DoubleLimitCovRiemannSums} (if one converges then the other does and vice-versa). Concerning expectations, using  \eqref{Eq:MeanSelfIntegralConvergence} it follows that the self-integrability of $K_{Z,M}$ is equivalent to the convergence of the expectations of the Riemann sums to a limit independent of the chosen Riemann system.

It follows that if $\int_{D}ZdM$ is uniquely-defined, from Lemma \ref{Lemma:L2Convergence}  the left side in \eqref{Eq:DoubleLimitCovRiemannSums} converges to $\Var(\int_{D}ZdM)$, hence the right side converges to the same value independently of both Riemann systems chosen, thus $K_{Z,M}^{(2)}$ is quasi-self-integrable. Conversely, if $K_{Z,M}^{(2)}$ is quasi-self-integrable then the double-limit at left side in \eqref{Eq:DoubleLimitCovRiemannSums} converges to the same value independently of the Riemann systems. Applying then Lemma \ref{Lemma:L2Convergence} point $\ref{It:L2Convergence1}$, one proves using equal Riemann systems $ \RiemannSystem{I}{x}{J} = \RiemannSystem{E}{y}{\tilde{J}}$ in equation \eqref{Eq:DoubleLimitCovRiemannSums} that the Riemann sums associated to such system converge in $\LpOmegaAP{2}$ to a random variable. The independence of the Riemann system follows hence when using two different Riemann systems in \eqref{Eq:DoubleLimitCovRiemannSums} and then applying Lemma \ref{Lemma:L2Convergence} point $\ref{It:L2COnvergence2}$, concluding thus that the corresponding Riemann sums converge to the same random variable. It follows hence that $\int_{D}ZdM$ is uniquely-defined as a limit in $\LpOmegaAP{2}$. $\blacksquare$

It is worth mention that Theorem \ref{Theo:IntZMGaussian} does not necessarily hold in a non-Gaussian context. Consider the following counter-example: $D = [0,1)$, $Z(t) = B^{2}(t)$, with $B$ a Gaussian Brownian motion and $M = B'$ the distributional derivative of $B$, which is a Gaussian White Noise. Then, we can use the explicit link between the Itô Integral $\int_{[0,1)}B^{2}(t)dB(t)$ and the Stratonovich integral $\int_{[0,1)}B^{2}(t)\circ dB(t) $:
\begin{equation}
\label{Eq:DifItoStratonovich}
\int_{[0,1)}B^{2}(t)\circ dB(t)  - \int_{[0,1)}B^{2}(t)dB(t) = \int_{[0,1)}B(t)dB(t),
\end{equation}
where we have used the traditional notation for Itô and Stratonovich integrals \citep[Equation 5.5]{kloeden1995numerical}. The Itô integral $\int_{[0,1)}B(t)dB(t)$ is not null, hence there is a difference between both integrals, thus $\int_{D}ZdM$ is not uniquely-defined. However, in this case $K_{B^{2},B'} = 0 $ (impair moments of a centred Gaussian vector are involved, which are always null), which is trivially self-integral, as well as $K_{B^{2},B'}^{(2)}$.

\subsection{Stochastic and self-integration over subsets}
\label{Sec:Subsets}

Both definitions of the self-integral and the uniquely-defined stochastic integral we have used are focused on a fixed integration domain $D$. Form the definitions, it is not clear this implies the possibility of integrating over subsets of $D$. In this Section we tackle such problem.

\begin{Theo}
\label{Theo:StochasticIntegralSubSet}
Let $D \in \BoundedBorelRd$. If the stochastic integral $\int_{D}ZdM$ is uniquely-defined in the sense of $L^{p}\OmegaAP$ for some $p \in [1,\infty]$, then for every $A \in \Borel{D}$ the stochastic integral $\int_{A}ZdM$ is uniquely-defined in the sense of $L^{p}\OmegaAP$ and the additivity property holds:
\begin{equation}
\label{Eq:StochasticIntegralAdditivity}
\int_{A}ZdM + \int_{D\setminus A} Z dM = \int_{D} Z dM.
\end{equation}
\end{Theo}

We shall use the following Lemma on sequences and sub-sequences on Banach spaces.

\begin{Lemma}
\label{Lemma:SequencesBanach}
Let $(E , \| \cdot \|_{E})$ be a Banach space. Let $(a_{n})_{n \in \mathbb{N}}$ and $(b_{n})_{n \in \mathbb{N}}$ be two sequences in $E$ satisfying the following property: there exists a unique $\ell \in E$ such that for every growing function $\sigma : \mathbb{N} \to \mathbb{N}$,
\begin{equation}
\label{Eq:PropertySequencesBanach}
\lim_{n \to \infty} a_{n} + b_{\sigma(n)} = \ell.
\end{equation}
Then, both $(a_{n})_{n \in \mathbb{N}}$ and $(b_{n})_{n \in \mathbb{N}}$ converge.
\end{Lemma}

\textbf{Proof of Lemma \ref{Lemma:SequencesBanach}:} For any two growing functions $\sigma , \sigma' : \mathbb{N} \to \mathbb{N}$, one has from property \eqref{Eq:PropertySequencesBanach}
\begin{equation}
\label{Eq:bnbmGoTo0}
\lim_{n \to \infty} b_{\sigma(n)} - b_{\sigma'(n)} = \lim_{n \to \infty} a_{n}+ b_{\sigma(n)} - (a_{n} + b_{\sigma'(n)} ) = \ell - \ell = 0.
\end{equation}
Suppose now that $(b_{n})_{n \in \mathbb{N}}$ does not converge. Since $E$ is Banach, $(b_{n})_{n}$ must not be a Cauchy sequence. There exists thus $\epsilon > 0 $ such that for every $n \in \mathbb{N}$ one can find $\tilde{n}, \tilde{m} \geq n$ such that $\| b_{\tilde{n}} - b_{\tilde{m}} \|_{E} \geq \epsilon$. We can hence construct two growing functions $\sigma_{1} , \sigma_{2} : \mathbb{N} \to \mathbb{N}$ as following. For $n = 0 $, we take $\sigma_{1}(n), \sigma_{2}(n) \geq 1 $ such that $\| b_{\sigma_{1}(0)} - b_{\sigma_{2}(0)} \|_{E} \geq \epsilon$. Then, for every $n \geq 1 $ we take $\sigma_{1}(n) , \sigma_{2}(n) > \max\lbrace \sigma_{1}(n-1) , \sigma_{2}(n-1) \rbrace$ such that $\| b_{ \sigma_{1}(n) } - b_{\sigma_{2}(n)} \|_{E} \geq \epsilon $. It follows that both $\sigma_{1} $ and $\sigma_{2}$ are growing and that $\| b_{\sigma_{1}(n)} - b_{\sigma_{2}(n)} \|_{E} \geq \epsilon$ for all $n \in \mathbb{N}$. But this contradicts the result \eqref{Eq:bnbmGoTo0}. $(b_{n})_{n \in \mathbb{N}}$ must hence be Cauchy, and thus convergent since $E$ in Banach. From property \eqref{Eq:PropertySequencesBanach} it follows immediately that $(a_{n})_{n \in \mathbb{N}}$ must also converge. $\blacksquare$

\textbf{Proof of Theorem \ref{Theo:StochasticIntegralSubSet}:} Let $A \in \Borel{D}$. Let $\RiemannSystem{A}{a}{J^{A}}$ be any Riemann system of $A$, and let $\RiemannSystem{\tilde{A}}{\tilde{a}}{J^{\tilde{A}}}$ be any Riemann system of $D\setminus A$. Then, one can construct a new Riemann system of $D$, say $\RiemannPartitionTag{I}{x}{J}$ through,
\begin{equation}
\label{Eq:SumOfRiemannSystems}
J_{n} = J_{n}^{A} \cup J_{n}^{\tilde{A}}; \quad  I_{j}^{n} =\begin{cases}
A_{j}^{n} & \hbox{ if } j \in J_{n}^{A} \\
\tilde{A}_{j}^{n} & \hbox{ if } j \in J_{n}^{\tilde{A}}
\end{cases}; \quad x_{j}^{n} = \begin{cases}
a_{j}^{n} & \hbox{ if } j \in J_{n}^{A} \\
\tilde{a}_{j}^{n} & \hbox{ if } j \in J_{n}^{\tilde{A}}
\end{cases}.
\end{equation} 
The Riemann sum associated to such system is the sum of the Riemann sums over $A$ and $D\setminus A$:
\begin{equation}
\label{Eq:RiemannSumsSumOfRiemannSystems}
\sum_{j \in J_{n}} Z(x_{j}^{n})M(I_{j}^{n}) = \sum_{j \in J_{n}^{A}} Z(a_{j}^{n})M(A_{j}^{n}) + \sum_{j \in J_{n}^{\tilde{A}}} Z(\tilde{a}_{j}^{n})M(\tilde{A}_{j}^{n}).  
\end{equation}
If $\sigma : \mathbb{N} \to \mathbb{N}$ is a growing function, then $\left( (\tilde{A}_{j}^{\sigma(n)})_{j \in J_{\sigma(n)}^{\tilde{A}}}  , (\tilde{a}_{j}^{\sigma(n)})_{j \in J_{\sigma(n)}^{\tilde{A}}}  \right)_{n\in \mathbb{N}} $ is a sub-Riemann system of $D\setminus A$. Using it together with $\RiemannSystem{A}{a}{J^{A}}$ for constructing a new Riemann system of $D$ analogously as in \eqref{Eq:RiemannSumsSumOfRiemannSystems}, since the stochastic integral $\int_{D} Z dM$ is uniquely-defined we must have 
\begin{equation}
\lim_{n \to \infty} \sum_{j \in J_{n}^{A}} Z(a_{j}^{n})M(A_{j}^{n}) + \sum_{j \in J_{\sigma(n)}^{\tilde{A}}} Z(\tilde{a}_{j}^{\sigma(n)})M(\tilde{A}_{j}^{\sigma(n)}) = \int_{D}ZdM.  
\end{equation}
From Lemma \ref{Lemma:SequencesBanach} using the Banach space $\LpOmegaAP{p}$, the Riemann sum $\sum_{j \in J_{n}^{A}} Z(a_{j}^{n})M(A_{j}^{n}) $ must converge. Set $\ell_{A} := \lim_{n \to \infty} \sum_{j \in J_{n}^{A}} Z(a_{j}^{n})M(A_{j}^{n})$. If $\RiemannSystem{\tilde{A}}{\tilde{a}}{J^{\tilde{A}}}$ denotes now any arbitrary Riemann system of $D\setminus A$, then by constructing a Riemann system of $D$ as in \eqref{Eq:SumOfRiemannSystems}, it follows
\begin{equation}
\label{Eq:LimitsAdditivityStochastic}
\lim_{n \to \infty}\sum_{j \in J_{n}^{\tilde{A}}} Z(\tilde{a}_{j}^{n})M(\tilde{A}_{j}^{n}) = \lim_{n \to \infty} \sum_{j \in J_{n}} Z(x_{j}^{n})M(I_{j}^{n}) - \sum_{j \in J_{n}^{A}} Z(a_{j}^{n})M(A_{j}^{n}) = \int_{D}Zd_{M} - \ell_{A}.
\end{equation}
Since the Riemann system of $D\setminus A$ was arbitrary, it follows that the stochastic integral $\int_{D\setminus A} Z dM$ is uniquely-defined and equals $\int_{D}Zd_{M} - \ell_{A}$. By taking arbitrary Riemann systems of $A$ and $D\setminus A$ and constructing an associated Riemann system of $D$ through \eqref{Eq:SumOfRiemannSystems}, it follows that the stochastic integral $\int_{A}ZdM$ is also uniquely-defined and equals $\ell_{A}$. The additivity property \eqref{Eq:StochasticIntegralAdditivity} follows from \eqref{Eq:LimitsAdditivityStochastic}. $\blacksquare$

The next Theorem is a \textit{résumé} of the main results exposed until now and it includes thus self-integrability over subsets for cross-covariance kernels. It also makes precise the important relations of expectation and covariances between stochastic integrals.

\begin{Theo}
\label{Theo:Resume}
Suppose $Z$ and $M$ are jointly Gaussian and centred. Let $D \in \BoundedBorelRd$. Then, the following  statements are equivalent:
\begin{enumerate}[(i)]
\item $K_{Z,M}$ is self-integrable over $D$ and $K_{Z,M}^{(2)}$ is quasi-self-integrable over $D\times D$.
\item $K_{Z,M}$ is self-integrable over every $A \in \Borel{D}$ and $K_{Z,M}^{(2)}$ is quasi-self-integrable over every set of the form $A\times B $ with $A,B \in \Borel{D}$.
\item The stochastic integral $\int_{D}ZdM$ is uniquely-defined as a limit in $\LtwoOmega$.
\item For every $A \in \Borel{D}$, the stochastic integral $\int_{A}ZdM$ is uniquely-defined as a limit in $\LtwoOmega$. 
\end{enumerate}
Moreover, if any of these statements holds, the expectation and covariance formulae hold:
\begin{equation}
\label{Eq:MeanStochasticIntegralA}
\mathbb{E}\left( \int_{A}ZdM \right) = \int_{A}K(x,dx), \quad \forall A \in \Borel{D},
\end{equation}
\begin{equation}
\label{Eq:CovStochasticIntegral}
\Cov\left( \int_{A}ZdM  , \int_{B}ZdM \right) = \int_{A\times B} C_{Z} dC_{M} + \tilde{\int}_{A\times B} K^{(2)}(x,dx), \quad \forall A , B \in \Borel{D}.
\end{equation}
\end{Theo}

\textbf{Proof of Theorem \ref{Theo:Resume}:} $(i) \Leftrightarrow (iii)$ is Theorem \ref{Theo:IntZMGaussian}. $(iii) \Leftrightarrow (iv)$ is Theorem \ref{Theo:StochasticIntegralSubSet}. $(ii) \Rightarrow (i)$ is obvious. We need only to prove $(iv) \Rightarrow (ii)$ and the formulae \eqref{Eq:MeanStochasticIntegralA} and \eqref{Eq:CovStochasticIntegral}. The self-integrability of $K_{Z,M}$ over $A$ and formula \eqref{Eq:MeanStochasticIntegralA} come from Proposition \ref{Prop:SelfIntegralNecessary}. For $K^{(2)}_{Z,M}$, let $\RiemannSystem{A}{a}{J^{A}}$ and $\RiemannSystem{B}{b}{J^{B}}$ be two arbitrary Riemann systems of $A$ and $B$  respectively. Since the associated Riemann sums approaching the stochastic integrals converge to $\int_{A}ZdM$ and $\int_{B}ZdM$ respectively, we use the same arguments as in \eqref{Eq:CovDoubleRiemannSum} and we conclude
\small
\begin{equation}
\Cov\Big( \int_{A}ZdM , \int_{B}ZdM \Big) = \int_{A\times B} C_{Z} dC_{M} + \lim_{(n,m) \to \infty} \sum_{j \in J_{n}^{A}}\sum_{k \in J_{m}^{B}} K_{Z,M}^{(2)}( (a_{j}^{n} , b_{k}^{m} ) , A_{j}^{n}\times B_{k}^{m} ).
\end{equation}
\normalsize
The double-limit at the right side holds then  for any arbitrary Riemann systems of $A$ and $B$. It follows that $K_{Z,M}^{(2)}$ is quasi-self-integrable over $A\times B$, and its quasi-self-integral satisfies equation \eqref{Eq:CovStochasticIntegral}. $\blacksquare$

Let us now prove that a self-integrable kernel $K$ is also self-integrable on subsets, without asking cross-pos-def properties or any other condition on $K^{(2)}$.

\begin{Theo}
\label{Theo:KSelfIntSubsets}
Let $K$ be an arbitrary cf-m kernel over $\Rd$, and let $D \in \BoundedBorelRd$. If $K$ is self-integrable over $D$, then it is self-integrable over every subset $A \in \Borel{D}$ and the additivity property holds:
\begin{equation}
\label{Eq:AdditivitySelfIntegral}
\int_{A}K(x,dx) + \int_{D\setminus A} K(x,dx) = \int_{D}K(x,dx), \quad \forall A \in \Borel{D}.
\end{equation}
\end{Theo}

\textbf{Proof of Theorem \ref{Theo:KSelfIntSubsets}:} The arguments are the same as in Theorem \ref{Theo:StochasticIntegralSubSet}, so we give the outline and leave details to the reader: use an arbitrary Riemann system of $A$ and an arbitrary Riemann system of $D \setminus A$. Construct a Riemann system of $D$ through \eqref{Eq:SumOfRiemannSystems}. The sum of the respective Riemann sums converges to $\int_{D}K(x,dx)$, as well as the sums when taking a sub-sequence of the Riemann sum of $D\setminus A$. Apply Lemma \ref{Lemma:SequencesBanach} in the case $E = \R$ to prove that the Riemann sums over $A$ and $D\setminus A$ converge. Argue then by varying the Riemann system of $D\setminus A$, that $K$ is self-integrable over $D\setminus A$, and conclude then that $K$ is self-integrable over $A$, with the additivity property \eqref{Eq:AdditivitySelfIntegral} holding. $\blacksquare$

For $K^{(2)}$, we only mention that, if $K$ is cross-pos-def, self-integrable over $D$ and such that $K^{(2)}$ is quasi-self-integrable over $D\times D$, then $K^{(2)}$ is quasi-self-integrable over every $A\times B$, $A,B \in \Borel{D}$. This comes directly from Theorem \ref{Theo:Resume}. In this case we have an additivity property for the quasi-self-integrals
\small
\begin{equation}
\tilde{\int}_{A\times A}K^{(2)}(x,dx) + \tilde{\int}_{A\times (D\setminus A)}K^{(2)}(x,dx) + \tilde{\int}_{(D\setminus A)\times A}K^{(2)}(x,dx) + \tilde{\int}_{(D\setminus A)\times(D\setminus A)}K^{(2)}(x,dx) = \tilde{\int}_{D\times D}K^{(2)}(x,dx),
\end{equation}
\normalsize
which can be obtained through equation \eqref{Eq:CovStochasticIntegral}, using $\int_{D}ZdM = \int_{A}ZdM + \int_{D\setminus A}ZdM$ and eliminating correspondingly the integrals of the form $\int C_{Z} dC_{M}$ that there appear.

\subsection{Measure structure}
\label{Sec:MeasureStructure}

Results in Section \ref{Sec:Subsets} have allowed us to define, under suitable self-integrability conditions over $D$, the applications $A \in \Borel{D} \mapsto \int_{A}K(x,dx)$, and $(A, B) \in \Borel{D}\times\Borel{D} \mapsto \tilde{\int}_{A\times B}K^{(2)}(x,dx)$. As we have seen, the first application is additive and it is not difficult to prove that the second is additive in each component when the other is fixed. However, it is not clear if $\sigma$-additivity properties can be concluded for those applications, that is, if the first one defines a measure over $D$ or if the second defines a measure over $D\times D$ (or a bi-measure whatsoever). This question is not without any use, since in such case, one can construct stochastic integrals which determine random measures. In particular, the following Corollary follows immediately from equations \eqref{Eq:MeanStochasticIntegralA} and \eqref{Eq:CovStochasticIntegral}.

\begin{Corol}
\label{Corol:IntZMRandomMeasure}
Suppose $Z$ and $M$ are jointly Gaussian. Then, the application $A  \mapsto \int_{A}ZdM$ determines a uniquely-defined random measure over $D$ if and only if the self-integrals $A \mapsto \int_{A}K_{Z,M}(x,dx)$ and the quasi-self-integrals $(A,B) \mapsto \tilde{\int}_{A\times B} K_{Z,M}^{(2)}(x,dx)$ are well-defined and determine measures over $D$ and $D\times D$ respectively.
\end{Corol}

It is then important to determine cases where a measure structure is induced in the self-integration process. Let us begin with the case of the self-integral of $K$. Rather than proving a $\sigma$-additivity property, we will construct directly the integrals of measurable functions with respect to the potential measure. Suppose $D$ is compact and $K$ is self-integrable over $D$. If $\varphi \in \mathcal{M}_{B}(D)$ is a simple function, say of the form $ \varphi = \sum_{i} a_{i} \mathbbm{1}_{A_{i}} $, then we can define using Theorem \ref{Theo:KSelfIntSubsets},
\begin{equation}
\label{Eq:DefIntSelfIntSimple}
\int_{D} \varphi(x)K(x,dx) := \sum_{i \in I} a_{i} \int_{A_{i}}K(x,dx).
\end{equation} 
In order to extent this to a more general $\varphi$, we start, as in the case of classical measures, with the case where $K$ is a \textit{positive} kernel over $D$, that is such that $K(x,A) \geq 0 $, for every $(x,A) \in D \times \Borel{D}$. in such case, one has for $\varphi$ simple
\begin{equation}
\label{Eq:SelfIntKcontinuousLinear}
\left| \int_{D}\varphi(x)K(x,dx) \right|  \leq \sum_{i} |a_{i}|\int_{A_{i}}K(x,dx) \leq \|\varphi \|_{\infty}\int_{D}K(x,dx). 
\end{equation}
For $\varphi \in C(D)$, we express it as a uniform limit of simple functions $(\varphi_{n})_{n}$ and we define $\int_{D}\varphi(x)K(x,dx) = \lim_{n \to \infty} \int_{D}\varphi_{n}(x)K(x,dx)$. From \eqref{Eq:SelfIntKcontinuousLinear} it follows that the limit is well-defined and it does not depend upon the sequence approaching $\varphi$. Equation \eqref{Eq:SelfIntKcontinuousLinear} also holds thus for $\varphi \in C(D)$, and the application $\varphi \in C(D) \mapsto \int_{D}\varphi(x)K(x,dx)$ defines a linear and continuous application in the supremum norm. From Riesz Representation Theorem (Remark \ref{Rem:RieszRepresentationD}), there exists a unique measure representing such application, and hence its definition can be extended to any $\varphi \in \mathcal{M}_{B}(D)$. We conclude the next Proposition.

\begin{Prop}
\label{Prop:KpositiveSelfIntMeasure}
If $K \geq 0$ is self-integrable over $D$ compact, then the application $A \mapsto \int_{A}K(x,dx)$ defines a (positive) measure. 
\end{Prop}

For the general case when $K$ is not positive, we consider an intuitive sufficient condition regarding the self-integrability of the total-variation kernel $|K|$. We remark that Definition \ref{Def:SelfIntegral} of the self-integral only applies for a cf-m kernel, hence we cannot in general apply it to $|K|$ since, as far as we now, $|K|$ may not be continuous. We need hence to add this extra requirement. Note that this anyway holds for $K$ cross-pos-def (Theorem \ref{Theo:KcrossPosDefEquicontinuous}).

\begin{Prop}
\label{Prop:KtotVarSelfIntegrable}
Suppose $D$ is compact. If $K$ is a cf-m kernel such that $|K|$ is also a cf-m kernel and if both $K$ and $|K|$ are self-integrable over $D$, then $A \mapsto \int_{A}K(x,dx)$ defines a measure satisfying
\begin{equation}
\label{Eq:SelfIntKlinearFunctionalKtotVar}
\left| \int_{D}\varphi(x)K(x,dx) \right| \leq \| \varphi \|_{\infty} \int_{D}|K|(x,dx), \quad \forall \varphi \in \mathcal{M}_{B}(D).
\end{equation}
\end{Prop}

\textbf{Proof:} If $\RiemannPartitionTag{A}{a}{J^{A}}$ is any Riemann system of $A$, then 
\begin{equation}
\label{Eq:ProvingSelfIntKAbouned}
\left|\int_{A}K(x,dx)\right| = \lim_{n \to \infty} \left| \sum_{j \in J_{n}^{A}} K(a_{j}^{n} , A_{j}^{n} ) \right| \leq  \lim_{n \to \infty}  \sum_{j \in J_{n}^{A}}  |K|(a_{j}^{n} , A_{j}^{n})  = \int_{A}|K|(x,dx). 
\end{equation}
It follows that for a simple function $\varphi = \sum_{i} a_{i} \mathbbm{1}_{A_{i}}$, one has
\begin{equation}
\label{Eq:ProvingSelfIntKvarphiContinuousKtotvar}
\left|\int_{D}\varphi(x)K(x,dx)\right| = \left| \sum_{i} a_{i} \int_{A_{i}}K(x,dx) \right| \leq \sum_{i}|a_{i}| \int_{A_{i}}|K|(x,dx) \leq \| \varphi \|_{\infty} \int_{D}|K|(x,dx).
\end{equation}
The linear application $\varphi \mapsto \int_{D}\varphi(x)K(x,dx) $ is hence continuous in the sense of the supremum norm. It determines hence a measure and its application can be extended to any $\varphi \in \mathcal{M}_{B}(D)$. $\blacksquare$

The general question of whether the self-integrability of $K$ implies or is implied by the self-integrability of $|K|$ remains open. Now, we can give similar results concerning the quasi-self-integral of $K^{(2)}$.

\begin{Prop}
\label{Prop:KpositiveK2quasiSelfIntMeasure}
Let $D \subset \Rd$ compact. If $K$ is positive, cross-pos-def and self-integrable over $D$, and if such that $K^{(2)}$ is quasi-self-integrable over $D\times D$, then $(A , B) \mapsto \int_{A\times B}K^{(2)}(x,dx)$ defines a measure over $D\times D$. 
\end{Prop}

\textbf{Proof:} We proceed similarly as in Proposition \ref{Prop:KpositiveSelfIntMeasure}. Of course $K$ positive implies $K^{(2)}$ positive. From Theorem \ref{Theo:Resume}, $K^{(2)}$ is quasi-self-integrable over every subset $A\times B$, $A,B \in \Borel{D}$. Hence, we can define a linear application over every simple function of the form $\psi = \sum_{i} \sum_{j} a_{i,j} \mathbbm{1}_{A_{i}}\otimes \mathbbm{1}_{B_{k}} $ as
\begin{equation}
\tilde{\int}_{D\times D} \psi(x)K^{(2)}(x,dx) := \sum_{i \in I} \sum_{j \in J} a_{i,j} \tilde{\int}_{A_{i}\times B_{j}} K(x,dx).
\end{equation}
And one can easily prove also following the same idea as in \eqref{Eq:SelfIntKcontinuousLinear} that for such $\psi$
\begin{equation}
\left| \tilde{\int}_{D\times D} \psi(x)K^{(2)}(x,dx) \right| \leq \| \psi \|_{\infty} \tilde{\int}_{D\times D} K^{(2)}(x,dx).
\end{equation}
Since every function in $C(D\times D)$ can be uniformly approached by functions of the form of $\psi$, one can apply the same approximation arguments and Riesz Representation Theorem to prove that the application $\psi \mapsto \tilde{\int}_{D\times D} \psi(x)K^{(2)}(x,dx)$ defines a measure and can be extended uniquely to any $\psi \in \mathcal{M}_{B}(D\times D)$. $\blacksquare$

The next Proposition is concluded following similar arguments which we left to the reader.

\begin{Prop}
\label{Prop:K2quasiSelfIntTotVar}
Suppose $D$ compact. Let $K$ be cross-pos-def, such that $K$ and $|K|$ are self-integrable over $D$ and $K^{(2)}$ and $|K^{(2)}|$ are quasi-self-integrable over $D\times D$. Then, the application $(A,B) \mapsto \tilde{\int}_{A\times B} K^{(2)}(x,dx)$ defines a measure over $D\times D$.
\end{Prop}

We finish this Section with the next Corollary which follows directly from Propositions \ref{Prop:KtotVarSelfIntegrable} and \ref{Prop:K2quasiSelfIntTotVar} and Theorem \ref{Theo:Resume}.

\begin{Corol}
\label{Corol:IntZMrandomMeasureTotVar}
Suppose $Z$ and $M$ jointly Gaussian and $D$ compact. If $K_{Z,M}$ and $|K_{Z,M}|$ are self-integrable over $D$, and if $K_{Z,M}^{(2)}$ and $|K_{Z,M}^{(2)}|$ are quasi-self-integrable over $D\times D$, then $A \mapsto \int_{A}ZdM$ defines a unique random measure over $D$ with mean and covariance measures given by \eqref{Eq:MeanStochasticIntegralA} and \eqref{Eq:CovStochasticIntegral}.
\end{Corol}

\begin{Rem}
\label{Rem:IntZMvarphiRandomMeasure}
When conditions in Corollaries \ref{Corol:IntZMRandomMeasure} or \ref{Corol:IntZMrandomMeasureTotVar} hold, since $\int_{A}ZdM$ defines a random measure, the integrals of every $\varphi \in \mathcal{M}_{B}(D)$ with respect to this measure are well-defined random variables, giving thus a unique definition for the stochastic integral
\begin{equation}
\label{Eq:IntZMvarphi}
\int_{D}Z(x)\varphi(x)dM(x),
\end{equation}
and such random variable should coincide with other intuitive forms of defining such integral. We remark that this allows, among many other things, to define properly distributional derivatives of the product $``ZM"$ when using $\varphi$ smooth with compact support.
\end{Rem}

\section{Applications}
\label{Sec:Applications}

In this Section we apply our results. We begin with two cases where conditions in Theorem \ref{Theo:IntZMGaussian} hold with a surprising convenience, allowing the unique-definition of some concepts involving stochastic integrals. We also present an important examples without uniquely-defined stochastic integrals, showing some limitations of our approach.

\subsection{Tensor product of Gaussian random measures}
\label{Sec:TensorProductOfGaussianRandomMeasure}

In Lemma \ref{Lemma:M1tensorM2indep} we have defined the tensor product of two independent random measures. Now, we will relax the hypothesis of independence but we will ask a Gaussianity condition. 
Let $M_{1}$ and $M_{2}$ be two real centred jointly Gaussian random measures over $\Rd$ and $\R^{p}$ respectively. Let us suppose, in addition, that the cross-covariance between them is determined by a measure $C_{M_{1} , M_{2}} \in \mathscr{M}(\Rd\times\R^{p})$, that is $\mathbb{E}(M_{1}(A)M_{2}(B) ) =  C_{M_{1} , M_{2}}(A \times B )$ for any $A \in \BoundedBorelRd, B \in \BoundedBorelRp$.  $M_{1}\otimes M_{2}$ should then be a new random measure over $\Rd\times\R^{p}$ satisfying $(M_{1}\otimes M_{2} )(A\times B ) = M_{1}(A)M_{2}(B)$. 

We follow the same approach as in the proof of Lemma \ref{Lemma:M1tensorM2indep}. Let $\psi \in C_{c}(\Rd\times\R^{p})$, with $\supp(\psi) \subset D_{d}\times D_{p}$, $D_{d}$, and $D_{p}$ compacts. We propose to define the action of $M_{1}\otimes M_{2}$ over $\psi$ through an iterated stochastic integral as in \eqref{Eq:DefM1tensorM2IteratedInt}. In order to verify that such stochastic integral is well-defined, we consider that $ Z_{M_{2}}(x):= \int_{\R^{p}} \psi(x,y) dM_{2}(y)$ is a mean-square continuous process over $\Rd$, null outside $D_{d}$. The cross-covariance kernel between $Z_{M_{2}}$ and $M_{1}$ is
\begin{equation}
\label{Eq:KZM2M1}
K_{Z_{M_{2}} , M_{1}}(x , A ) = \int_{A \times \R^{p}} \psi(x,v) dC_{M_{1},M_{2}}(u,v) = \int_{D_{d}\times D_{p}} \psi(x,v)\mathbbm{1}_{A}(u)dC_{M_{1},M_{2}}(u,v).
\end{equation}
If $\RiemannPartitionTag{I}{x}{J}$ is a Riemann system of $D_{d}$, the associated Riemann sum approaching the potential self-integral of $K_{Z_{M_{2}} , M_{1}}$ is
\begin{equation}
\label{Eq:RiemannSumSelfIntKZM2M1}
\sum_{j \in J_{n}} \int_{D_{d}\times D_{p}} \psi(x_{j}^{n},v)\mathbbm{1}_{I_{j}^{n}}(u)dC_{M_{1},M_{2}}(u,v) = \int_{D_{d}\times D_{p}} \left( \sum_{j \in J_{n}} \psi(x_{j}^{n},v)\mathbbm{1}_{I_{j}^{n}}(u) \right)dC_{M_{1},M_{2}}(u,v).
\end{equation}
From the uniform continuity of $\psi$ over $D_{d}\times D_{p}$, one has $\sum_{j \in J_{n}} \psi(x_{j}^{n},v)\mathbbm{1}_{I_{j}^{n}}(u) \to \psi(u,v) $ uniformly over $D_{d}\times D_{p}$. Since $C_{M_{1},M_{2}}$ is a measure, this implies the convergence of the integrals to a limit independent of the Riemann system used. $K_{Z_{M_{2}} , M_{1}}$ is hence self-integrable over $D_{d}$ with self-integral
\begin{equation}
\label{Eq:SelfIntKZM2M1}
\int_{D_{d}}K_{Z_{M_{2}} , M_{1}}(x,dx) = \int_{\Rd\times \R^{p}} \psi(x,y) dC_{M_{1},M_{2}}(x,y).
\end{equation}
The second order kernel is given by (computations left to the reader)
\begin{equation}
\label{Eq:K2ZM2M1}
K_{Z_{M_{2}},M_{1}}^{(2)}( (x,y) , E ) =\int_{D_{d}\times D_{p}}\int_{D_{d}\times D_{p}} \psi(x,\tilde{v})\mathbbm{1}_{E}(u,\tilde{u}) \psi(y,v) d C_{M_{1},M_{2}}(u,v) d C_{M_{1},M_{2}}(\tilde{u} , \tilde{v} ).
\end{equation}
$K_{Z_{M_{2}},M_{1}}^{(2)}$ is actually self-integrable over $D_{d}\times D_{d}$ in the sense of Definition \ref{Def:SelfIntegral}, which implies of course quasi-self-integrability. To see this, we argue similarly as in the case of the self-integrability of $K_{Z_{M_{2}} , M_{1}}$  using the uniform continuity of $\psi$ which allows to conclude that for any Riemann system $\RiemannSystem{E}{e}{J}$ of $D_{d}\times D_{d}$, (say $e_{j}^{n} = ( x_{j}^{n} , y_{j}^{n} ) \in D_{d}\times D_{d}$) one has 
\begin{equation}
\label{Eq:UnifConvApproachingCrossKernelK2ZM2M1}
\sum_{j \in J_{n}} \psi(x_{j}^{n} , \tilde{v}) \mathbbm{1}_{E_{j}^{n}}(u,\tilde{u})\psi(y_{j}^{n} , v) \xrightarrow[n \to \infty]{} \psi(u,\tilde{v})\psi(\tilde{u} , v)
\end{equation}
uniformly in $(u,v,\tilde{u},\tilde{v}) \in D_{d}\times D_{p}\times D_{d}\times D_{p}$. The self-integral of $K_{Z_{M_{2}},M_{1}}^{(2)}$ is hence given by
\begin{equation}
\label{Eq:SelfIntegralK2ZM2M1}
\int_{D_{d}\times D_{d}}K_{Z_{M_{2}},M_{1}}^{(2)}(x,dx) = \int_{D_{d}\times D_{p} \times D_{d} \times D_{p}}  \psi(u,\tilde{v})\psi(\tilde{u} , v) d(C_{M_{1},M_{2}} \otimes  C_{M_{1},M_{2}} )( (u,v) , (\tilde{u} , \tilde{v} ) ).
\end{equation}
Theorem \ref{Theo:IntZMGaussian} can be hence applied. The stochastic integrals of the form
\begin{equation}
\label{Eq:M1tensroM2PsiGeneral}
\langle M_{1}\otimes M_{2} , \psi \rangle := \int_{\Rd}\int_{\R^{p}}\psi(x,y) dM_{2}(y) dM_{1}(x)
\end{equation}
for $\psi \in C_{c}(\Rd\times\R^{p})$ are uniquely-defined. Let us describe the mean and covariance of such variables. For the mean, it follows from \eqref{Eq:SelfIntKZM2M1} that
\begin{equation}
\label{Eq:MeanM1tensorM2}
\mathbb{E}\left( \langle M_{1} \otimes M_{2} , \psi \rangle \right) = \int_{\Rd\times\R^{p}} \psi(x,y) dC_{M_{1},M_{2}}(x,y).
\end{equation}
For the covariance consider for simplicity $\psi_{1}$ and $\psi_{2}$ real. We then compute through Riemann sums and using Lemma \ref{Lemma:Cov4MomentGaussian}
\small
\begin{equation}
\label{Eq:CovM1tensorM2}
\begin{aligned}
\Cov&\left( \langle M_{1} \otimes M_{2} , \psi_{1} \rangle , \langle M_{1} \otimes M_{2} , \psi_{2} \rangle  \right) \\
&=  \lim_{(n,m)\to \infty} \sum_{j \in J_{n}}\sum_{k \in \tilde{J}_{m}} \Cov\left( \int_{\R^{p}}\psi_{1}(x_{j}^{n} , v)dM_{2}(v)M_{1}(I_{j}^{n}) ,  \int_{\R^{p}}\psi_{2}(y_{k}^{m} , \tilde{v})dM_{2}(\tilde{v})M_{1}(E_{k}^{m}) \right) \\
&= \lim_{(n,m)\to \infty} \sum_{j \in J_{n}}\sum_{k \in \tilde{J}_{m}} \int_{\R^{p}\times\R^{p}} \psi_{1}(x_{j}^{n},v)\psi_{2}( y_{k}^{m} , \tilde{v} ) dC_{M_{2}}(v,\tilde{v}) C_{M_{1}}(I_{j}^{n}\times E_{k}^{m}) \\
&\quad \quad +\sum_{j \in J_{n}}\sum_{k \in \tilde{J}_{m}} \int_{\Rd\times \R^{p}} \psi_{1}(x_{j}^{n} , v) \mathbbm{1}_{E_{k}^{m}}(u) dC_{M_{1},M_{2}}(u,v)\int_{\Rd\times \R^{p}} \psi_{2}(y_{k}^{m} , v) \mathbbm{1}_{I_{j}^{n}}(\tilde{u}) dC_{M_{1},M_{2}}(\tilde{u},\tilde{v}) \\
&= \int_{\Rd\times\Rd\times\R^{p}\times\R^{p}} \psi_{1}(\tilde{u},v)\psi_{2}(u,\tilde{v}) d(C_{M_{1}}\otimes C_{M_{2}} )( u,\tilde{u} , v , \tilde{v} ) \\
&\quad \quad + \int_{\Rd\times\R^{p}\times\Rd\times\R^{p}} \psi_{1}(\tilde{u},v)\psi_{2}(u,\tilde{v}) d(C_{M_{1},M_{2}}\otimes C_{M_{1},M_{2}} )( u, v , \tilde{u}, \tilde{v} ).
\end{aligned}
\end{equation}
\normalsize
The last double-limit is justified using similar uniform-approximation arguments as in \eqref{Eq:UnifConvApproachingCrossKernelK2ZM2M1}. Since $C_{M_{1}}, C_{M_{2}}$ and $C_{M_{1},M_{2}}$ are all measures, it follows that expression \eqref{Eq:CovM1tensorM2} defines a measure over $\Rd\times\R^{p}\times\Rd\times\R^{p}$. We can express these results in a more synthetic manner as
\begin{equation}
\label{Eq:MeanCovM1tensorM2}
m_{M_{1}\otimes M_{2}} := C_{M_{1},M_{2}} \quad ; \quad C_{M_{1}\otimes M_{2}} = C_{M_{1}}``\otimes" C_{M_{2}} + C_{M_{1},M_{2}} ``\otimes" C_{M_{1},M_{2}},
\end{equation}
where $C_{M_{1}}``\otimes" C_{M_{2}}$  and $C_{M_{1},M_{2}} ``\otimes" C_{M_{1},M_{2}}$ denote respectively the only measures over $\Rd\times\R^{p} \times \Rd \times \R^{p}$ such that for all $A_{1},A_{2} \in \BoundedBorelRd, B_{1},B_{2} \in \mathcal{B}_{B}(\R^{p})$ one has
\begin{equation}
( C_{M_{1}}``\otimes" C_{M_{2}} )(A_{1}\times B_{1} \times A_{2} \times B_{2} ) = C_{M_{1}}(A_{1}\times A_{2})C_{M_{2}}(B_{1}\times B_{2});
\end{equation}
\begin{equation}
(C_{M_{1},M_{2}} ``\otimes" C_{M_{1},M_{2}})( A_{1}\times B_{1} \times A_{2} \times B_{2}  ) = C_{M_{1},M_{2}}(A_{1}\times B_{2}) C_{M_{1},M_{2}}(A_{2}\times B_{1}).
\end{equation}
Note that the measure $ C_{M_{1}}``\otimes" C_{M_{2}} $ is not the measure  $C_{M_{1}} \otimes C_{M_{2}} $ since the latter is defined over $\RdxRd\times\R^{p}\times\R^{p}$. As well, the measure $C_{M_{1},M_{2}} ``\otimes" C_{M_{1},M_{2}}$ is not $C_{M_{1},M_{2}} \otimes C_{M_{1},M_{2}}$ since the former has the cross-product property of connecting $A_{1}$ with $B_{2}$ rather than with $B_{1}$. Using Theorem \ref{Theo:ExtensionRandomMeasure} and thus extending the definition of $M_{1}\otimes M_{2}$ to functions in $\mathcal{M}_{B,c}(\Rd\times\R^{p})$,  we conclude the following Theorem.
\begin{Theo}
\label{Theo:M1tensorM2}
Let $M_{1}$ and $M_{2}$ be jointly Gaussian random measures over $\Rd$ and $\R^{p}$ respectively. Suppose that the cross-covariance structure between $M_{1}$ and $M_{2}$ is determined by a measure $C_{M_{1},M_{2}} \in \mathscr{M}(\Rd\times \R^{p})$. Then, there exists a unique random measure over $\Rd\times\R^{p}$, noted $M_{1}\otimes M_{2}$, which satisfies
\begin{equation}
(M_{1}\otimes M_{2})(A\times B) = M(A)M(B), \quad \forall A \in \BoundedBorelRd, B \in \BoundedBorel{\R^{p}}.
\end{equation}
\end{Theo}

Consider the iterated integrals in \eqref{Eq:M1tensroM2PsiGeneral} in the other sense $\int_{\R^{p}} \int_{\Rd} \psi(x,y)dM_{1}(x) dM_{2}(y)$. Following the same arguments as before, one can conclude that such iterated integral is uniquely-defined. The mean and the covariance structures of the application $\psi \mapsto \int_{\R^{p}} \int_{\Rd} \psi(x,y)dM_{1}(x) dM_{2}(y)$ can be computed analogously as we did for \eqref{Eq:MeanCovM1tensorM2}, obtaining that both structures are determined by measures. In addition, it is clear that $\int_{\R^{p}} \int_{\Rd} \psi(x,y)dM_{1}(x) dM_{2}(y) $ coincide with $\int_{\R^{d}} \int_{\R^{p}} \psi(x,y)dM_{2}(y) dM_{1}(x)$ for every $\psi \in C_{c}(\Rd) \otimes C_{c}(\R^{p})$. Since $C_{c}(\Rd) \otimes C_{c}(\R^{p})$ is dense in $C_{c}(\Rd\times\R^{p})$ with the uniform-on-compacts topology, and since both applications define continuous linear functions from $C_{c}(\Rd\times\R^{p})$ to $\LtwoOmega$, both must coincide over the whole space $C_{c}(\Rd\times\R^{p})$. We conclude the following Fubini Theorem.

\begin{Theo}[\textbf{Fubini Theorem for random measures}]
\label{Theo:FubiniTensor}
Let $M_{1}$ and $M_{2}$ be jointly Gaussian random measures over $\Rd$ and $\R^{p}$ respectively such that the cross-covariance structure between $M_{1}$ and $M_{2}$ is determined by a measure $C_{M_{1},M_{2}} \in \mathscr{M}(\Rd\times \R^{p})$. Then, for every $\psi \in C_{c}(\Rd\times\R^{p})$ the triple equality
\begin{equation}
\int_{\R^{p}} \int_{\Rd} \psi(x,y)dM_{1}(x) dM_{2}(y) = \langle M_{1}\otimes M_{2} , \psi \rangle = \int_{\Rd}\int_{\R^{p}}\psi(x,y) dM_{2}(y) dM_{1}(x)
\end{equation}
holds, the iterated integrals being uniquely-defined.
\end{Theo}

It is worth mentioning that Theorem \ref{Theo:FubiniTensor} works only for $\psi$ continuous with compact support. While it also works (trivially) for functions in $\mathcal{M}_{B,c}(\Rd)\otimes \mathcal{M}_{B,c}(\R^{p})$, \textit{it does not work in general for any $\psi \in \mathscr{M}_{B,c}(\Rd\times\R^{p})$}. We show a counter-example in Section \ref{Sec:DerivativeRandomMeasures}. The main issue is that in such case the iterated integrals may not be uniquely-defined, even if the broad tensor product measure is. One way of expressing this fact is the following: \textit{tensor products of Gaussian random measures do not have a canonical stochastic disintegration}.

We finish by mentioning that through the use of tensor products one can define the convolution between two Gaussian random measures, at least in the case where both $M_{1}$ and $M_{2}$ are finite over $\Rd$ ($C_{M_{1}}$ and $C_{M_{2}}$ have finite total mass over $\Rd\times\Rd$) and a measure $C_{M_{1},M_{2}}$ describes their cross-covariance. In such case, the application $\varphi \in C_{c}(\Rd) \mapsto \langle M_{1}\otimes M_{2} , (x,y) \mapsto \varphi(x+y) \rangle$ should define a new random measure. This is an advance in the research of a convolutive algebra between random measures, at least in the Gaussian case, which may contribute to the questions proposed in \citep[Chapter 5]{rao2012random} and \citep{park2016random}.

\subsection{Fractional Brownian motion and its derivative}
\label{Sec:FracBrownianMotion}

A fractional Brownian motion over $[0,1]$ with Hurst index $H \in (0,1)$ is a centred Gaussian process $(B_{H}(t))_{t \in [0,1]}$ with covariance function
\begin{equation}
\label{Eq:CovFractionalBrowMotion}
C_{B_{H}}(t,s) = \frac{|t|^{2H} + |s|^{2H} - |t-s|^{2H}}{2}.
\end{equation}
Suppose $H > \frac{1}{2}$. If $B'_{H}$ is the distributional derivative of $B_{H}$, then its covariance is \citep{borisov2006constructing}
\begin{equation}
\label{Eq:CovDerFracBrownian}
C_{B_{H}'} = \frac{\partial^{2}C_{B_{H}}}{\partial t \partial s}  = H(2H-1)|t-s|^{2H-2}.
\end{equation}
Since $\frac{1}{2} < H < 1 $, we have $-1 < 2H - 2 < 0 $. It follows that covariance \eqref{Eq:CovDerFracBrownian} is identified with a function over $[0,1] \times [0,1]$ which is integrable but not continuous. $B_{H}'$ is hence a random measure with absolutely continuous covariance measure 
\begin{equation}
\label{Eq:CovDerFracBrownianMeasure}
C_{B_{H}'}(E) = H(2H-1)\int_{E} \frac{d(u,v)}{|u-v|^{2-2H}} , \quad \forall E \in \Borel{[0,1]\times[0,1]}.
\end{equation} 
Some authors also call $B_{H}'$ a long-range dependence process, see \citep{gay1990class,anh1999possible}. Now, using $B_{H}(t) = B_{H}'([0,t])$ we can compute the cross-covariance kernel between $B_{H}$ and its derivative $B_{H}'$:
\begin{equation}
\label{Eq:CrossCovBhBhder}
\begin{aligned}
K_{B_{H},B_{H}'}(t,A) = \mathbb{E}\left( B_{H}'([0,t])B_{H}'(A) \right) &= H(2H-1)\int_{A\times[0,t]} \frac{d(u,v)}{|u-v|^{2-2H}} \\
&= H\int_{A} u^{2H-1} + |t-u|^{2H-1} \sgn(t-u) du. 
\end{aligned}
\end{equation}
Let $\RiemannSystem{I}{t}{J}$ be any Riemann system of $[0,1]$. Then,
\begin{equation}
\label{Eq:RiemannSumSelfIntKBhBhder}
\sum_{j \in J_{n}} K_{B_{H},B_{H}'}(t_{j}^{n},I_{j}^{n}) = H \int_{0}^{1} u^{2H-1}du + H\sum_{j \in J_{n}} \int_{I_{j}^{n}} |t_{j}^{n} - u|^{2H-1} \sgn(t_{j}^{n} - u ) du. 
\end{equation}
For the second term at the right side of \eqref{Eq:RiemannSumSelfIntKBhBhder}, we consider that $2H-1 > 0 $ and hence
\begin{equation}
\begin{aligned}
\left| \sum_{j \in J_{n}} \int_{I_{j}^{n}} |t_{j}^{n} - u|^{2H-1} \sgn(t_{j}^{n} - u ) du \right| &\leq \sum_{j \in J_{n}} \diam(I_{j}^{n})^{2H-1} \lambda(I_{j}^{n}) \\
&\leq \max_{j \in J_{n}} \diam(I_{j}^{n})^{2H-1}  \xrightarrow[n \to \infty]{} 0.
\end{aligned}
\end{equation}
It follows that $K_{B_{H},B_{H}'}$ is self-integrable with self-integral
\begin{equation}
\label{Eq:SelfIntegralKBhBhder}
\int_{[0,1]\times[0,1]}K_{B_{H},B_{H}'}(x,dx) = H \int_{0}^{1} u^{2H-1}du = \frac{1}{2}.
\end{equation}
On the other hand, the second order kernel $K_{B_{H},B_{H}'}^{(2)}$ is given by
\small
\begin{equation}
\label{Eq:K2BhBhder}
K_{B_{H},B_{H}'}^{(2)}((t,s) , A\times B) = H^{2}\int_{A\times B} \left( u^{2H-1} + |s-u|^{2H-1}\sgn(s-u) \right)\left( v^{2H-1} + |t-v|^{2H-1}\sgn(t-v)\right) d(u,v).
\end{equation}
\normalsize
If $\RiemannSystem{D}{s}{\tilde{J}}$ is another Riemann system of $[0,1]$, the double-Riemann sums approaching the quasi-self-integral of $K_{B_{H},B_{H}'}^{(2)}$ are given by the integral over $[0,1]\times[0,1]$ of the double-sequence of functions
\small
\begin{equation}
\label{Eq:MultiSeqApproachingK2BhBhder}
(u,v) \mapsto H^{2}\sum_{j \in J_{n}}\sum_{k\in \tilde{J}_{m}}\mathbbm{1}_{I_{j}^{n}}(u) \mathbbm{1}_{D_{k}^{m}}(v)\left( u^{2H-1} + |s_{k}^{m}-u|^{2H-1}\sgn(s_{k}^{m}-u) \right)\left( v^{2H-1} + |t_{j}^{n}-v|^{2H-1}\sgn(t_{j}^{n}-v)\right).
\end{equation}
\normalsize
Let us study the convergence of the functions \eqref{Eq:MultiSeqApproachingK2BhBhder} when $(n,m) \to \infty$. Fix $u \in [0,1]$ and for each $n$ denote $j_{u} \in J_{n}$ the only element of $J_{n}$ such that $u \in I_{j_{u}}^{n}$. It follows that $x_{j_{u}}^{n} \to u $ when $n \to \infty$. We do analogously with any $v \in [0,1]$ and $k_{v} \in \tilde{J}_{m}$. Thus \eqref{Eq:MultiSeqApproachingK2BhBhder} can be re-written as
\begin{equation}
\label{Eq:MultiSeqApprK2BhBhderIndexuv}
(u,v) \mapsto H^{2}\left( u^{2H-1} + |s_{k_{v}}^{m}-u|^{2H-1}\sgn(s_{k_{v}}^{m}-u) \right)\left( v^{2H-1} + |t_{j_{u}}^{n}-v|^{2H-1}\sgn(t_{j_{u}}^{n}-v)\right).
\end{equation}
If $u \neq v$, since the $\sgn$ function is continuous outside $0$, expression \eqref{Eq:MultiSeqApprK2BhBhderIndexuv} converges as $(n,m) \to \infty$ to
\begin{equation}
\label{Eq:LimitFunctionApproachK2BhBhder}
H^{2}\left( u^{2H-1} + |v-u|^{2H-1}\sgn(v-u) \right)\left( v^{2H-1} + |u-v|^{2H-1}\sgn(u-v)\right).
\end{equation}
This last function is integrable over $[0,1]\times[0,1]$. In addition, the functions \eqref{Eq:MultiSeqApprK2BhBhderIndexuv} are all bounded by the constant $H^{2}(1 + \sqrt{2}^{2H-1})^{2} $, hence they converge almost-everywhere and dominated to the function \eqref{Eq:LimitFunctionApproachK2BhBhder}. It follows from LDCT Theorem for double-sequences (see Lemma \ref{Lemma:LDCTdoubleSequence}), that the double-sequence of integrals of \eqref{Eq:MultiSeqApproachingK2BhBhder} over $[0,1]\times[0,1]$ converges to
\begin{equation}
\label{Eq:ConvergenceToQuasiSelfIntK2BhBhder}
H^{2}\int_{0}^{1}\int_{0}^{1}\left( u^{2H-1} + |v-u|^{2H-1}\sgn(v-u) \right)\left( v^{2H-1} + |u-v|^{2H-1}\sgn(u-v)\right)d(u,v). 
\end{equation}  
This integral, which is then the quasi-self-integral of $K_{B_{H},B_{H}'}^{(2)}$, can be computed using a long but elementary calculus. The final result is
\begin{equation}
\label{Eq:QuasiSelfIntegralK2BhBhder}
\tilde{\int}_{[0,1]\times[0,1]}K_{B_{H},B_{H}'}^{(2)}(x,dx) = \frac{1}{2} -H\mathbf{B}(2H , 2H+1) - \frac{H}{2(4H-1)},
\end{equation}
where $\mathbf{B}$ is the beta function. From Theorems \ref{Theo:IntZMGaussian} and \ref{Theo:Resume}, the stochastic integrals
\begin{equation}
\label{Eq:StochInfBhBhder}
\int_{A}B_{H}(t)dB_{H}'(t), \quad A \in \Borel{[0,1]}
\end{equation}
are uniquely-defined. In addition, $K_{B_{H},B_{H}'}$ is positive, so we are in the quite regular case of Corollary \ref{Corol:IntZMrandomMeasureTotVar}. The application $A \mapsto \int_{A}B_{H}(t)dB_{H}'(t)$ defines then a random measure over $[0,1]$.

This property of fractional Brownian motion with $H > \frac{1}{2}$ makes a huge contrast with the case of Browninan motion $(H=\frac{1}{2})$ briefly presented in the proof of Proposition \ref{Prop:ExampleNotDefSelfInt}. At least in what concerns the definition of simple stochastic integrals, the stochastic calculus with respect to $B_{H}$ is better behaved than with respect to $B$.

\subsection{Processes whose derivatives are random measures}
\label{Sec:DerivativeRandomMeasures}

Let $Z$ and $M$ be such that $M$ is the (distributional) derivative of $Z$. The example studied in Section \ref{Sec:FracBrownianMotion} is a particular case, and other important case is the one of Brownian motion and its derivative White Noise described partially in the proof of Proposition \ref{Prop:ExampleNotDefSelfInt}. Let $D = [0,1]$. We can start from $M$ and define
\begin{equation}
\label{Eq:ZprimitiveM}
Z(t):= M([0,t]),
\end{equation}
which is the only (distributional) primitive of $M$ such that $Z(0) = M(\lbrace 0 \rbrace)$. In general, such $Z$ is mean-square càdlàg (right-continuous with left limits). Let us suppose, for simplicity, that $M$ is such that $Z$ is mean-square continuous, hence $M$ must not have atoms ($M(\lbrace t \rbrace) \stackrel{a.s.}{=} 0 $ for every $t \in [0,1]$). The cross-covariance kernel is 
\begin{equation}
\label{Eq:KZprimitiveM}
K_{Z,M}(t,A) = C_{M}( [0,t]\times A ).
\end{equation}
In this case, Riemann sums approaching the self-integral are always bounded:
\begin{equation}
\label{Eq:KRiemannSumSelfIntBoundedZprimM}
\left| \sum_{j \in J_{n}} C_{M}((0, x_{j}^{n}] , I_{j}^{n})  \right| \leq \sum_{j \in J_{n}} |C_{M}|(  D \times I_{j}^{n} ) = |C_{M}|(D\times D).
\end{equation}
By taking a sub-sequence, one can always find a Riemann sum which converges, but in general this does not imply self-integrability. 

Let us consider an important particular case which is widely used in theory in applications: $M$ is an \textit{orthogonal random measure}. By this, we mean that $M$ has a covariance of the form 
\begin{equation}
\label{Eq:CovOrthogonal}
C_{M}(A\times B ) = \nu_{M}(A \cap B), \quad \forall A,B \in \BoundedBorelRd,
\end{equation}
being $\nu_{M}$ a positive measure over $\Rd$. Orthogonal random measures take non-correlated values over disjoint subsets. Some authors restrain the study of random measure to this class, sometimes with the stronger requirement of independence in disjoint subsets \citep{kingman1967completely,passeggeri2020extension}. Sadly, in this case the stochastic integral $\int_{[0,1]}ZdM$ is not uniquely-defined. Indeed, in this case the cross-covariance kernel \eqref{Eq:KZprimitiveM} is given by
$ K_{Z,M}(t,A) = \nu_{M}([0,t]\cap A)$, and the same problem as in Proposition \ref{Prop:ExampleNotDefSelfInt} is present. Consider partitions of $[0,1]$ consistent of left-closed intervals $(I_{j}^{n})_{j \in J_{n}}$, with $a_{j}^{n} := \inf I_{j}^{n} $ and $b_{j}^{n} := \sup I_{j}^{n}$. Note that since $M$ has no atoms, $\nu_{M}$ has no atoms as well. If we take as tag-points the initial points $x_{j}^{n} = a_{j}^{n}$, one has
\begin{equation}
\sum_{j \in J_{n}} K_{Z,M}(a_{j}^{n} , I_{j}^{n} ) = \sum_{j \in J_{n}} \nu_{M}([0,a_{j}^{n}]\cap I_{j}^{n} ) =\sum_{j \in J_{n}} \nu_{M}(\lbrace a_{j}^{n}\rbrace ) = 0.
\end{equation}
But using the uniform continuity of the function $ t \mapsto \nu_{M}( [0,t] ) $ over $[0,1]$, one can always find, for every $n$, tag-points $(x_{j}^{n})_{j \in J_{n}}$ such that $\sum_{j \in J_{n}} \nu_{M}( I_{j}^{n}\cap (x_{j}^{n} , b_{j}^{n} ] ) \to 0$ as $n \to \infty$. Hence,
\begin{equation}
\sum_{j \in J_{n}} K_{Z,M}(x_{j}^{n} , I_{j}^{n} ) = \sum_{j \in J_{n}} \nu_{M}([a_{j}^{n} , x_{j}^{n} ]) = \sum_{j \in J_{n}} \nu_{M}( I_{j}^{n} ) - \nu_{M}( I_{j}^{n}\cap (x_{j}^{n} , b_{j}^{n} ]  ) \xrightarrow[n \to \infty]{} \nu_{M}(D).
\end{equation}
It follows that kernels of the type \eqref{Eq:KZprimitiveM} are never self-integrable when $M$ is orthogonal without atoms. This last result shows one of the many difficulties of developing and stochastic calculus with respect to important random measures. While this may look disappointing for some readers, there is a curious analysis we can do when $M$ is a Gaussian random measure which explains this phenomenon. Since $M\otimes M$ is a well-defined random measure(Theorem \ref{Theo:M1tensorM2}, $C_{M}$ being the cross-covariance measure in this case), the integral $\langle M\otimes M , \psi \rangle$ is uniquely-defined for every $\psi \in \mathcal{M}_{B}([0,1]\times[0,1])$. This includes $\psi(t,s) =  \mathbbm{1}_{[0,t]}(s)$. One could be thus tempted to say
\begin{equation}
``\int_{[0,1]} M([0,t]) dM(t) = \int_{[0,1]\times [0,1]} \mathbbm{1}_{[0,t]}(s)  d(M \otimes M)(s,t)".
\end{equation}
But integral at the left side is not uniquely-defined. In addition, since $M$ is orthogonal one always has
\begin{equation}
\int_{[0,1]\times[0,1]}  \mathbbm{1}_{[0,t]}(s) d(M\otimes M)(s,t) \neq \int_{[0,1]\times[0,1]} \mathbbm{1}_{[0,t)}(s)d(M\otimes M)(s,t). 
\end{equation}
Indeed, since the mean measure of $M\otimes M $ is $C_{M}$ and $C_{M}$ has the form \eqref{Eq:CovOrthogonal}, we have
\begin{equation}
\mathbb{E}\left( \int_{[0,1]\times[0,1]}  \mathbbm{1}_{[0,t]}(s) d(M\otimes M)(s,t) \right) = \int_{[0,1]} \mathbbm{1}_{[0,t]}(t) d\nu_{M}(t) = \nu_{M}([0,1]),
\end{equation}
while
\begin{equation}
\mathbb{E}\left(   \int_{[0,1]\times[0,1]} \mathbbm{1}_{[0,t)}(s)d(M\otimes M)(s,t) \right) = \int_{[0,1]} \mathbbm{1}_{[0,t)}(t) d\nu_{M}(t) = 0.
\end{equation}
This implies that if we try to interpret the integrals over $[0,1]\times[0,1]$ as iterated integrals, we obtain
\begin{equation}
\label{Eq:IteratedIntegralsM}
\int_{[0,1]} M([0,t]) dM(t) \neq \int_{[0,1]} M([0,t)) dM(t),
\end{equation} 
which is not so intuitive since the processes $t \mapsto M([0,t])$ and $t \mapsto M([0,t))$ are essentially the same. While $\mathbbm{1}_{[0,t]}(s)$ and $\mathbbm{1}_{[0,t)}(s)$ are almost-everywhere equal, this holds with respect to the Lebesgue measure over $[0,1]\times[0,1]$, but the measure that matters here is $C_{M}$. The latter is concentrated on the diagonal $\lbrace t = s \rbrace$, thus with respect to $C_{M}$ the functions $\mathbbm{1}_{[0,t]}(s)$ and $\mathbbm{1}_{[0,t)}(s)$ are entirely different.

It follows that well-definiteness of the tensor measure $M\otimes M$ over measurable bounded functions does not imply the unique-definiteness of the iterated integrals \eqref{Eq:IteratedIntegralsM}. This is a counter-example to Fubini Theorem \ref{Theo:FubiniTensor}, showing that it statement cannot be extended to any $\psi$. This peculiar result seems more clear when we realize that \textit{Riemann sums, while they are an intuitive way of approaching the integrals in \eqref{Eq:IteratedIntegralsM}, they are a poor manner of approaching integral $\langle M\otimes M , \psi \rangle$}. Indeed, when we use Riemann sums approaching integrals \eqref{Eq:IteratedIntegralsM}, what we are actually doing is
\begin{equation}
\sum_{j \in J_{n}} M([0,x_{j}^{n}])M(I_{j}^{n}) = \int_{[0,1]\times[0,1]} \sum_{j \in J_{n}} \mathbbm{1}_{[0,x_{j}^{n}]}(t)\mathbbm{1}_{I_{j}^{n}}(s)  d( M\otimes M )(t,s).
\end{equation}
The sequence of functions $\sum_{j \in J_{n}} \mathbbm{1}_{[0,x_{j}^{n}]}\otimes \mathbbm{1}_{I_{j}^{n}}$ is such that for $(t,s) \in [0,1]\times[0,1]$, if $j_{s} \in J_{n}$ is the index such that $s \in I_{j}^{n}$, then $\sum_{j \in J_{n}} \mathbbm{1}_{[0,x_{j}^{n}]}(t) \mathbbm{1}_{I_{j}^{n}}(s) = \mathbbm{1}_{[0,x_{j_{s}}^{n}]}(t)$. Since $x_{j_{s}}^{n} \to s $ as $ n \to \infty$, if $t > s$  then  $\mathbbm{1}_{[0,x_{j_{s}}^{n}]}(t) \to 1$. If $t < s$, then $\mathbbm{1}_{[0,x_{j_{s}}^{n}]}(t) \to 0 $. But when $t=s$ (where $C_{M}$ is concentrated), the limit depends upon the tag-points $(x_{j}^{n})_{j \in J_{n}}$. If for some $s$, $(x_{j_{s}}^{n})_{n}$ is a sequence which has sub-sequences both at the right and left side of $s$, then $\mathbbm{1}_{[0,x_{j_{s}}^{n}]}(t)$ does not converge. Is this lack of clearness and universality of the convergence of functions of the form $\sum_{j \in J_{n}} \mathbbm{1}_{[0,x_{j}^{n}]}\otimes \mathbbm{1}_{I_{j}^{n}}$ what actually impeaches us to uniquely-define stochastic integrals the way we have tried.

\section{Perspectives and Comments}
\label{Sec:Perspectives}

We present some open questions and future possible applications of the results here exposed. They concern essentially conditions of self-integrability of kernels and the construction of products of more general stochastic structures.

\subsection{Some open questions}
\label{Sec:OpenQuestions}

Maybe one of the first questions some readers will have in mind is weather we can find some practical conditions to verify self-integrability of a cf-m kernel. In this work we focused on the connection between the self-integral and uniquely-defined stochastic integrals, but we have not focused on criteria for self-integrability. Besides the case of kernels in $C(D)\otimes \mathscr{M}(D)$ and in its completion with the projective topology which we have mentioned quickly in the introduction (without details), we have not given a general criterion for determining if a kernel is self-integrable or not. In a forthcoming paper we will describe some of those in terms of Karhunen-Loève expansions of both $Z$ and $M$, see \citep{carrizov2022karhunen}.

One notion which is also interesting is whether for any given kernel there exists a Riemann system whose associated Riemann sum converges. In  such case, \textit{a} self-integral can be obtained regardless if the use of another Riemann system would provide the same value. Thus, the question is we could find a \textit{Riemann-system-selection convention} which would always allow the convergence of the Riemann sums. It is not clear if for any arbitrary kernel $K$ one can construct such Riemann system. There are some cases where the existence of a converging Riemann sum is guaranteed, for instance using boundedness arguments as in equation \eqref{Eq:KRiemannSumSelfIntBoundedZprimM}, but this does not apply in general. For instance, in Appendix \ref{Sec:AppendixUnbounded} we show an example of a cross-pos-def kernel over $[-1,1]$ for which the use of Riemann systems using same-length intervals form always unbounded Riemann sums. Note that this question is not without use, since it is related to the possibility of defining a stochastic integral, even if non-uniquely, for a wide variety of integrands and integrators (Itô and Stratonovich constructions are examples Riemann-system-selection conventions).

The self-integrability of the second order kernel $K^{(2)}$ is also a property to be analysed. For instance, does the self-integrability of $K$ implies or is implied by the quasi-self-integrability of $K^{(2)}$?. We have not found a counter-example where one is self-integrable and the other is not. The same question stands for the self-integrability of the total-variation kernel $|K|$.

Another question concerns the measure structure of the studied objects. In Section \ref{Sec:MeasureStructure} we have remarked that a measure structure for the applications $A \mapsto \int_{A}K(x,dx)$ and $A \mapsto \int_{A}ZdM$ are not very simple to justify, giving only strong sufficient conditions for it even in the cases where self-integrability and unique-definiteness are present. Counter-examples have not been found. This question is crucial for knowing we can apply to the self-integral and the stochastic integral the same operation we know for measures and random measures, helping in the construction of a potential stochastic calculus grounded on these concepts.

Finally we consider a tensor-product approach to stochastic integration as a possible application of Theorem \ref{Theo:M1tensorM2}. Over $[0,1]\subset \R$, if $Z$ and $M$ are jointly Gaussian and $Z$ is a process whose derivative is a random measure $M_{Z}$, then the integral $\int_{[0,1]\times [0,1]}\mathbbm{1}_{[0,t]}(s) d(M_{Z}\otimes M)(s,t)$ is uniquely-defined provided the existence of the cross-covariance measure $C_{M_{Z},M}$. If so, integrals of the form $t \mapsto \int_{[0,t]} Z(u)dM(u) $ can be defined through a precise convention: the tensor product measure $M_{Z}\otimes M$ is unique, there is no ambiguity in such definition. The resulting process will likely have a mean-square càdlàg behaviour. A stochastic calculus in distributional sense can be hence developed, and it will have the advantage of being regular with respect to approximations in its tensor-product form: if $\mathbbm{1}_{[0,t]}(s)$ is approached by continuous functions, for instance, then the stochastic integrals of such functions will converge to the final stochastic integral. This regularity can be thus exploited to develop calculus and approximation rules.

\subsection{Generalized stochastic tensor products and products}
\label{Sec:GeneralizedTensorProducts}

This work may open a door for exploring the unique-definition of products between more abstract stochastic structures. If $\mu \in \mathscr{M}(\Rd)$ and $f \in C(\Rd)$, one can always define the product $f\mu$ as a new measure over $\Rd$ through $A \mapsto \int_{A}f d\mu $. As we have remarked, when $f$ and $\mu$ are stochastic their product can be canonically defined only if the corresponding stochastic integral is uniquely-defined. The same idea works in a more abstract context: if $T$ is a distribution over $\Rd$ ($T \in \DistributionsRd$) and $f \in C^{\infty}(\Rd)$, the product $fT$ is a new distribution given by $\varphi \in \DRd \mapsto \langle T , f\varphi \rangle$. Suppose now $T$ is a random distribution, also called a generalized random field (GeRF) \citep{gelfand1964generalized,ito1954stationary}, that is, a continuous linear operator $T: \DRd \to \LpOmega{2}$. Suppose $Z$ is a stochastic process with smooth covariance function. Then, given any $\varphi \in \DRd$, one could try to define $``\langle T , Z\varphi \rangle"$, but here we will have the same problem of non-uniqueness. Following an analogous approach of Riemann sums, one could try to approach $Z\varphi$ with a smooth process of the form $\sum_{j \in J_{n}} Z(x_{j}^{n})\varphi(x_{j}^{n})\phi_{j}^{n}$, with conveniently selected tag-points $x_{j}^{n}$ and functions $(\phi_{j}^{n})_{j \in J_{n}} \subset\DRd$, and then do $\langle T , Z\varphi \rangle:= \lim_{n \to \infty} \sum_{j \in J_{n}} Z(x_{j}^{n})\varphi(x_{j}^{n}) \langle T , \phi_{j}^{n}
\rangle$. In general, such limit will depend upon the selection of the tag-points and the functions $\phi_{j}^{n}$. Note that in this case we also have a cross-covariance kernel, $K_{Z\varphi,T} : \Rd\times\DRd \to \mathbb{R}$ given by $K_{Z\varphi,T}(x,\phi) = \Cov( Z(x)\varphi(x) , \langle T , \phi \rangle )$, which satisfies $K_{Z\varphi,T}( x , \cdot ) \in \DistributionsRd$ and $K_{Z\varphi,T}(\cdot , \phi ) \in \DRd$. Hence, the analogue of the self-integral \eqref{Eq:SelfIntegralDefRiemmanSums} would be the application of the right-component of $K_{Z\varphi,T}$, which is a distribution, to the left-component of $K_{Z\varphi,T}$, which is a smooth function with compact support. The possibility of this \textit{self-application}, which is also a \textit{trace operator}, will be hence linked with the possibility of defining uniquely the stochastic application $\langle T , Z\varphi \rangle$, in a completely analogous way to the link between self-integrals and uniquely-defined stochastic integrals in Theorem \ref{Theo:IntZMGaussian}. This approach can also be used when multiplying other stochastic structures with well-defined deterministic analogue: tempered distributions with smooth functions with polynomially bounded derivatives, primitives of measures with $C^{1}$ functions, etc.

One can also use this approach to define tensor products of Gaussian GeRFs, analogously to Section \ref{Sec:TensorProductOfGaussianRandomMeasure}. If $T,S$ are jointly Gaussian GeRFs over $\Rd$ and $\R^{p}$ respectively, the tensor product $T\otimes S$ should always be a well-defined GeRF over $\Rd\times\R^{p}$. If $\psi \in \mathscr{D}(\Rd\times\R^{p})$, the process $x \mapsto \langle S , \phi(x, \cdot) \rangle$ is mean-square smooth with compact support, and hence $\langle T \otimes S , \psi \rangle$ could be defined as the iterated application $\langle T ,x \mapsto \langle S , \phi(x, \cdot) \rangle$. This last product is expected to be uniquely-defined since, analogously to the criterion required in Theorem \ref{Theo:M1tensorM2}, the cross-covariance structure between two GeRFs is always identified to a distribution in $\mathscr{D}'(\Rd\times\R^{p})$ (Schwartz kernel's Theorem, see \citep[Chapter 51]{treves1967topological} and \citep[Theorem V.12]{reed1980methods}). Stochastic convolution can also be tackled through this approach.

It is then expected that this work will partially enlighten the analysis of stochastic products. Diverse sophisticated theories of constructing such products have been developed the last decades \citep{holden2009stochastic,gubinelli2004controlling,gubinelli2015paracontrolled,hairer2014theory}. The approach followed in this work is remarkable for its simplicity and elementaryness with respect to such theories. However, as the reader may already have thought, it can only be applied in some particular cases with unique-integrability properties. Let us say it: the most interesting applications (Itô and Stratonovich calculus, for instance), rely mainly upon the construction of stochastic integrals which are not uniquely-defined. Our approach may help hence in some circumstances but not always. Never mind, the big moral of this work, and which can always be taken into account when defining any kind of product or non-linear operation between two stochastic objects is the following: \textit{the impossibility of defining such products does not rely only on the regularities of the factors but also on their inter-dependence structure.}

\begin{appendices}

\section{Proof of Theorem \ref{Theo:ExtensionRandomMeasure}}
\label{Sec:AppendixProofTheoRieszRandom}

We will first prove the following Lemma which is a LDCT for double-sequences of functions.
\begin{Lemma}
\label{Lemma:LDCTdoubleSequence}
Let $(E,\mathscr{E}, \mu)$ be a measure space $(\mu \geq 0)$. Let $(f_{n,m})_{n,m \in \mathbb{N}}$ be a double-sequence of complex functions such that $\sup_{n,m \in \mathbb{N}} |f_{n,m}| \in \mathscr{L}^{1}(E,\mathscr{E},\mu)$. Suppose that the double limit $\lim_{n,m \to \infty} f_{n,m} $ exists $\mu$-almost everywhere. Then, the $\mu-$almost everywhere defined function $f(x) = \lim_{n,m \to \infty} f_{n,m}(x)$ is in $\mathscr{L}^{1}(E,\mathscr{E},\mu)$ and
\begin{equation}
\label{Eq:L1Convergefnmf}
\lim_{(n,m) \to \infty}\int_{E}|f_{n,m}-f|d\mu = 0.
\end{equation}
\end{Lemma}

\textbf{Proof of Lemma \ref{Lemma:LDCTdoubleSequence}:}\footnote{This proof has been borrowed from the StackExchange discussion \url{https://math.stackexchange.com/questions/448931/dominated-convergence-thm-dct-for-double-sequences}, consulted for the last time March the 6th 2023.} $f \in \mathscr{L}^{1}(E,\mathscr{E},\mu)$ is obvious since $\sup_{n,m \in \mathbb{N}} |f_{n,m}| \in \mathscr{L}^{1}(E,\mathscr{E},\mu)$. If we suppose that the double limit \eqref{Eq:L1Convergefnmf} does not converge to $0$, then there exists $\epsilon > 0 $ such that for any $k \in \mathbb{N}$ there exists $n_{k} , m_{k} \geq k$ such that
\begin{equation}
\label{Eq:SubseqIntBiggerEpsilon}
\int_{E}|f_{n_{k},m_{k}}-f|d\mu \geq \epsilon. 
\end{equation}
It follows that we can construct the sequence of functions $f_{k} := f_{n_{k},m_{k}}$, which convergence $\mu$-almost surely to $f$ and such that $\sup_{k \in \mathbb{N}}|f_{k}|$ is integrable. By the traditional LDCT \citep[VI.9]{doob1953stochastic}, we have $f_{k} \to f $ in $\mathscr{L}^{1}(E,\mathscr{E},\mu)$. But this contradicts \eqref{Eq:SubseqIntBiggerEpsilon}.  $\blacksquare$

We will also use the following Theorem of extension of vector valued measures.

\begin{Theo}
\label{Theo:ExtensionKupka}
Let $\mathcal{A}$ be an algebra of subsets of a set $D$, $\sigma(\mathcal{A})$ be the $\sigma$-algebra generated by $\mathcal{A}$ and let $X$ be a Banach space. Let $\mu : \mathcal{A} \to X $ be a $\sigma$-additive function on $\mathcal{A}$. Then, there exists a unique measure $\mu^{*} : \sigma(\mathcal{A}) \to X$  for which $\mu = \mu^{*}$ over $\mathcal{A}$, if and only if for every sequence of pairwise disjoint sets $(E_{n})_{n \in \mathbb{N}} \subset \mathcal{A}$ it holds $\lim_{n \to \infty} \mu(E_{n}) = 0 $.
\end{Theo}
\textbf{Proof of Theorem \ref{Theo:ExtensionKupka}:} See \citep{kupka1978caratheodory} for a clear exposition of this theorem and  \citep{sion1969outer} for the detailed proof. $\blacksquare$

\textbf{Proof of Theorem \ref{Theo:ExtensionRandomMeasure}:} From the linearity of the application $\varphi \in \mapsto \langle m_{M} , \varphi  \rangle$ and the bi-linearity of $(\varphi,\phi) \mapsto \langle C_{M} , \varphi \otimes \phi \rangle$, one proves easily the linearity of $M$. $M$ is also continuous since
\begin{equation}
\label{Eq:McontinuousC(D)L2}
\mathbb{E}(|\langle M , \varphi \rangle|^{2} ) =  |\langle m_{M} , \varphi \rangle|^{2} + \langle C_{M} , \varphi\otimes \overline{\varphi} \rangle \leq \left( |m_{M}|(D)^{2} + |C_{M}|(D\times D)\right) \| \varphi \|_{\infty}^{2}.
\end{equation}
Let us extend the definition of $M$ to a more general $\varphi$. Consider a uniformly bounded sequence $(\varphi_{n})_{n \in \mathbb{N}} \subset C(D)$ converging point-wise to some function $\varphi$. Then,
\begin{equation}
\label{Eq:MvarphinCauchyL2}
\mathbb{E}\left(|\langle M , \varphi_{n}-\varphi_{m} \rangle|^{2} \right) =  |\langle m_{M} , \varphi_{n} - \varphi_{m} \rangle |^{2} + \langle C_{M} , (\varphi_{n} - \varphi_{m})\otimes \overline{\varphi_{n} - \varphi_{m}} \rangle.
\end{equation}
Since $\lim_{n,m \to \infty} \varphi_{n}-\varphi_{m} = 0$ point-wise and dominated by a constant, and since both $m_{M}$ and $C_{M}$ are finite measures (since $D$ is compact), from Lemma \ref{Lemma:LDCTdoubleSequence} expression  \eqref{Eq:MvarphinCauchyL2} goes to $0$ as $n,m \to \infty$. $(\langle M , \varphi_{n} \rangle)_{n \in \mathbb{N}}$ is then a Cauchy sequence in $\LtwoOmega$ which converges to a unique random variable which we shall denote $\langle M , \varphi \rangle$. This limit does not depend upon the sequence $(\varphi_{n})_{n}$ chosen, since if we consider another one $(\phi_{n})_{n} \subset C(D)$ convergent to $\varphi$ point-wise and dominated, then LDCT allows to conclude that $\langle M , \varphi_{n} - \phi_{n} \rangle \to 0 $. It is also clear that the extension is linear in $\varphi$, and that the mean and covariance expressions \eqref{Eq:MeanCovIntegralM} hold. We remark now that any function of the form $\varphi = \mathbbm{1}_{R \cap D}$, where $R$ is a rectangle of $\Rd$, that is, a set of the form $R = I_{1}\times...\times I_{d}$, with $I_{1},... , I_{d}$ intervals of $\R$, can be approached by continuous functions in this way. It follows that we can define uniquely $\tilde{M}(R) := \langle M , \mathbbm{1}_{R \cap D} \rangle$ for every rectangle set $R$. Let us denote $\mathcal{R}(D)$ the algebra generated by the sets of the form $R\cap D$. From linearity, it is clear that $\tilde{M}$ is well defined over $\mathcal{R}(D)$, and since the mean and covariance structures are defined by measures, it is also $\sigma$-additive over $\mathcal{R}(D)$. Let $(E_{n})_{n \in \mathbb{N}} \subset \mathcal{R}(D)$ be an arbitrary sequence of mutually disjoint sets, and let $E := \bigcup_{n \in \mathbb{N}} E_{n} \in \Borel{D}$. Then,
\begin{equation}
\label{Eq:LimitMstronglyBounded}
\lim_{n \to \infty}\mathbb{E}(|\tilde{M}(E_{n})|^{2}) = m_{M}(E_{n})^{2} + C_{M}(E_{n}\times E_{n}) \leq |m_{M}|(E_{n})^{2} + |C_{M}|(D \times E_{n} ).
\end{equation}
Since both $m_{M}$ and $C_{M}$ are measures, we have
\begin{equation}
\label{Eq:SeriesmMBounded}
\sum_{n \in \mathbb{N}} |m_{M}|(E_{n}) = |m_{M}|(E) < \infty,
\end{equation}
\begin{equation}
\label{Eq:SeriesCmBounded}
\sum_{n \in \mathbb{N}} |C_{M}|(D\times E_{n}) = |C_{M}|(D\times E) < \infty.
\end{equation}
The series \eqref{Eq:SeriesmMBounded} and \eqref{Eq:SeriesCmBounded} are thus convergent, hence their terms must converge to $0$. By applying this in equation \eqref{Eq:LimitMstronglyBounded} we conclude $\lim_{n \to \infty} M(E_{n}) = 0 $ in the sense of the Banach space $\LtwoOmega$. Since $\Borel{D} = \sigma(\mathcal{R}(D))$, from Theorem \ref{Theo:ExtensionKupka}, there exists a unique extension of $\tilde{M}$ to $\Borel{D}$, defining a measure with values in $\LtwoOmega$.
The equality between $\langle M , \varphi \rangle$ and $\int_{D}\varphi d\tilde{M}$  is obtained from the construction of $\int_{D}\varphi d\tilde{M}$ using Riemann sums with Riemann systems whose partitions are formed by rectangles, which can always be done since any continuous function can be approached uniformly by linear combinations of indicators of rectangles. $\blacksquare$

\section{An unbounded Riemann sum}
\label{Sec:AppendixUnbounded}

In this section we consider a cross-pos-def kernel $K$ for which the Riemann sums approaching its self-integral using same-length intervals are necessarily unbounded. Hence, no possible self-integral defined through a partition with homogeneous intervals can be constructed. Consider the kernel defined over $D = [-1,1]$ by
\begin{equation}
\label{Eq:KernelUnboundedRS}
K(x,A) = \int_{A} \frac{du}{|x-u|^{\frac{1}{8}}}, \quad (x,A) \in [-1,1]\times\Borel{[-1,1]}.
\end{equation}
This kernel is the cross-covariance between a White Noise $W$ over $[-1,1]$ and the process
\begin{equation}
\label{Eq:DefZUnboundedRS}
Z(x) = \int_{-1}^{1} \frac{dW(u)}{|x-u|^{\frac{1}{8}}}.
\end{equation}
Note that for every $x \in [-1,1]$ the function $u \mapsto |x-u|^{-\frac{1}{8}}$ is in $L^{2}[-1,1]$, hence the integral with respect to the White Noise \eqref{Eq:DefZUnboundedRS} is well-defined. In addition, $Z$ is mean-square continuous since
\begin{equation}
\begin{aligned}
\mathbb{E}\left( [Z(x)-Z(x_{0}) ]^{2}  \right)  &= \mathbb{E}\left(  \left[ \int_{-1}^{1}  \frac{1}{|x-u|^{\frac{1}{8}}}  - \frac{1}{|x_{0}-u|^{\frac{1}{8}}}    dW(u) \right]^{2} \right) \\
&= \int_{-1}^{1} \left[  \frac{1}{|x-u|^{\frac{1}{8}}}  - \frac{1}{|x_{0}-u|^{\frac{1}{8}}}  \right]^{2} du \\
&= \int_{-1}^{1} \frac{\left[ |x\vee x_{0} - u |^{\frac{1}{8}} - |x \wedge x_{0} - u |^{\frac{1}{8}} \right]^{2}}{  |x \vee x_{0} - u |^{\frac{1}{4}}|x \wedge x_{0} - u |^{\frac{1}{4}}  } du.
\end{aligned}
\end{equation}
The numerator in last integral can be bounded by an arbitrarily small $\epsilon > 0 $ when $x \to x_{0}$ using the uniform continuity of the function $u \mapsto |u|^{\frac{1}{8}}$ over $[-1,1]$. For the denominator, the integral can be split in the intervals $[-1, \frac{x+x_{0}}{2} ]$ and $[ \frac{x+x_{0}}{2} , 1]$. In the first interval one has $|u - x\wedge x_{0}| \geq |u - x \vee x_{0} |$, so $|u - x\wedge x_{0}|^{-\frac{1}{4}} \leq |u - x \vee x_{0} |^{-\frac{1}{4}}$ and hence the integral over $[-1, \frac{x+x_{0}}{2} ]$ can be bounded by 
\begin{equation}
\epsilon \int_{-1}^{\frac{x+x_{0}}{2} } \frac{du}{|u-x\vee x_{0}|^{\frac{1}{2}}} = \epsilon 2 \left[ \left( \frac{|x-x_{0}|}{2}  \right)^{\frac{1}{2}} - (1 + x \vee x_{0} )^{\frac{1}{2}}  \right] \leq \epsilon 2(1+\sqrt{2}).
\end{equation}
An analogous analysis can be use to give an arbitrarily small bound to the integral over $[ \frac{x+x_{0}}{2} , 1]$. This proves the continuity of $Z$. Let $\RiemannPartitionTag{I}{x}{J}$ be a Riemann system of $[-1,1]$ for which for every $n$ the sets $I_{j}^{n}$ are intervals with same length $\ell_{n}$ (so $|J_{n}| = \frac{2}{\ell_{n}}$). Set $a_{j}^{n} := \inf(I_{j}^{n})$ and $b_{j}^{n} := \sup(I_{j}^{n})$. Then, the associated Riemann sum satisfies
\begin{equation}
\begin{aligned}
\sum_{j \in J_{n}} \int_{I_{j}^{n}} \frac{du}{|x_{j}^{n} - u |^{\frac{1}{8}}} &= \sum_{j \in J_{n}} \int_{a_{j}^{n}}^{x_{j}^{n}} \frac{du}{(x_{j}^{n} - u )^{\frac{1}{8}}} + \int_{x_{j}^{n}}^{b_{j}^{n}} \frac{du}{(u-x_{j}^{n} )^{\frac{1}{8}}} \\
&= \sum_{j \in J_{n}} \frac{8}{7}\left( (x_{j}^{n} - a_{j}^{n})^{\frac{7}{8}}  + (b_{j}^{n} - x_{j}^{n})^{\frac{7}{8}}   \right) \\
&\geq \sum_{j \in J_{n}}   \frac{8}{7} \left( \frac{\ell_{n}}{2} \right)^{\frac{7}{8}} = \frac{8}{7}|J_{n}|^{\frac{1}{8}} \nearrow \infty.
\end{aligned}
\end{equation}
This result proves that the Riemann sum is not bounded. Hence a self-integral using equal-length intervals cannot be constructed. The same holds for the construction of a stochastic integral between $Z$ and $W$. It is possible to prove, nonetheless, that when using non-equal-length intervals, one can find a convergent Riemann sums. Take for example a Riemann system where the $I_{j}^{n}$ are intervals, with $|J_{n}| = n $, and that every interval has length $e^{-n}$ except for one, which thus has length $2 - (n-1)e^{-n}$. Then the associated Riemann sum is bounded (details left to the reader), and hence it has a convergent sub-sequence.

\end{appendices}

\bibliography{mibib}
\bibliographystyle{apacite}

\end{document}